\documentclass[french,english]{amsart}

\usepackage{babel}
\usepackage[T1]{fontenc}
\usepackage{amsmath,amssymb}
\usepackage[all]{xy}
\usepackage{graphicx}

\newtheorem{propdef}{Proposition-D\'efinition}[section]
\newtheorem{lemme}[propdef]
{Lemme }
\newtheorem{theoreme}[propdef]
{Th\'eor\`eme }

\newtheorem{corollaire}[propdef]
{Corollaire }
\newtheorem{proposition}[propdef]{Proposition}

\theoremstyle{definition}
\newtheorem{ex}[propdef]{Exemple }
\newtheorem{rque}{Remarque}
\newtheorem{definition}[propdef]{D\'efinition }

\newcommand{\lgw}{\longrightarrow}
\newcommand{\lgm}{\longmapsto}
\newcommand{\ovl}{\overline}

\newcommand{\Frac}{\text{Frac}}

\newcommand{\ord}{\text{ord}}

\newcommand{\wdh}{\widehat}

\renewcommand{\l}{\lambda}
\renewcommand{\O}{\mathcal{O}}
\newcommand{\mfk}{\mathfrak}
\newcommand{\m}{\mathfrak{m}}
\newcommand{\I}{\mathfrak{I}}

\renewcommand{\k}{\Bbbk}

\newcommand{\pgcd}{\text{pgcd}}

\renewcommand{\mod}{\text{mod}}

\newcommand{\N}{\mathbb{N}}

\renewcommand{\a}{\alpha}
\renewcommand{\b}{\beta}

\newcommand{\e}{\varepsilon}

\begin{document}
\selectlanguage{french}
\title{Lemme d'Artin-Rees, th\'eor\`eme d'Izumi et fonction de Artin}

\author{Guillaume Rond}
\email{rond@picard.ups-tlse.fr}
\address{\newline
Laboratoire E. Picard - Universit\'e P. Sabatier\\
118, route de Narbonne, 
 31062 - Toulouse - Cedex 4 \\
 France}
\begin{abstract}
We interpret the Artin-Rees lemma and the Izumi theorem in term of Artin function and we obtain a stable version of the Artin-Rees lemma. We present different applications of these interpretations. First we show that the Artin function of $X_1X_2-X_3X_4$, as a polynomial in the ring of power series in more than three variables, is not bounded by an affine function. Then we prove that the Artin functions of a class of polynomials are bounded by affine functions and we use this to compute approximated integral closures of ideals.

\end{abstract}
\maketitle

\footnotetext[1]{\selectlanguage{english}2000 Mathematics Subject Classification: 13B40, (Primary), 13B22, 14B12 (Secondary)}
\selectlanguage{french}

\section{Introduction}
Nous rappelons quelques r\'esultats d'approximation, mais nous donnons tout d'abord la d\'efinition suivante :

\begin{definition}
Nous appellerons couple $(A,\,\mathfrak{I})$ la donn\'ee d'un anneau commutatif unitaire $A$ et d'un id\'eal $\I$ de $A$. Nous dirons que le couple $(A,\,\I)$ est n{\oe}th\'erien (resp. local, complet, r\'eduit, int\`egre) si l'anneau $A$ est n{\oe}th\'erien (resp. local, complet, r\'eduit, int\`egre).
\end{definition}
Nous pouvons alors d\'efinir les propri\'et\'es d'approximation et d'approximation forte :
\begin{definition}
Soit $(A,\,\mathfrak{I})$  un couple n{\oe}th\'erien et  $\widehat{A}$ le compl\'et\'e de $A$ pour la topologie $\I$-adique. Nous dirons que $(A,\,\mfk{I})$ v\'erifie \textit{la propri\'et\'e d'approximation} (PA) (resp. v\'erifie \textit{la propri\'et\'e d'approximation} pour $f$) si pour tout syst\`eme d'\'equations polynomiales not\'e $f(X)=0$ \`a coefficients dans $A$ (resp. si pour le syst\`eme d'\'equations polynomiales not\'e $f(X)=0$ \`a coefficients dans $A$), pour toute solution $\ovl{x}\in \wdh{A}$ et pour tout $i\in \mathbb{N}$, il existe une solution $x$ dans $A$ de ce syst\`eme qui v\'erifie $x_{j}=\ovl{x}_{j}\quad \mod\, \mfk{I}^{i+1} $ pour tout $j$.\\
Dans le cas o\`u $A$ est local et $\I$ est son id\'eal maximal, nous dirons que $A$ a la propri\'et\'e d'approximation.
\end{definition}

\begin{definition}
Soit $(A,\, \mfk{I})$ un couple n{\oe}th\'erien. Nous dirons que $(A,\, \mfk{I})$ v\'erifie \textit{la propri\'et\'e d'approximation forte} (PAF) si pour tout  syst\`eme d'\'equations polynomiales not\'e $f(X)=0$ \`a coefficients dans $A$, il existe une fonction \`a valeurs enti\`eres $\beta$ avec la propri\'et\'e suivante. Soient  $x\in A^{n}$ et $i\in\N$ tels que $$f(x)=0\quad \mod\, \mfk{I}^{\beta(i)+1}.$$
Alors il existe $\ovl{x}\in A^{n}$ tel que $$f(\ovl{x})=0\text{ et } x_j\equiv \ovl{x}_j\quad \mod \, \mfk{I}^{i+1}\,\text{pour tout } j.$$
La plus petite fonction vérifiant cette propri\'et\'e sera appel\'ee fonction de Artin de l'id\'eal $(f)$.\\
L\`a encore, si $A$ est local et $\I$ est son id\'eal maximal, nous dirons que $A$ a la propri\'et\'e d'approximation forte.
\end{definition}

\begin{rque}
Nous pouvons v\'erifier que les deux d\'efinitions pr\'ec\'edentes ne d\'ependent pas des g\'en\'erateurs de l'id\'eal $(f)$. Nous parlerons donc indiff\'eremment de syst\`eme d'\'equations polynomiales et d'id\'eal de $A[X]$.
\end{rque}

Nous avons les deux r\'esultats suivants :
\begin{theoreme}\cite{Ar}\cite{P-P}\cite{Spi2}\cite{Sw}
Soit $(A,\,\mfk{I})$ une paire hens\'elienne, n{\oe}th\'erienne. Alors $(A,\,\mfk{I})$ poss\`ede la propri\'et\'e d'approximation si  $A\lgw \widehat{A}$ est r\'egulier (o\`u $\wdh{A}$ est le compl\'et\'e $\mfk{I}$-adique de $A$).
\end{theoreme}
Nous rappelons qu'un morphisme d'anneaux n{\oe}th\'eriens $\varphi :A\lgw B$ est dit r\'egulier si il est plat et si pour tout id\'eal premier $P$ de $A$, la fibre $B\otimes_A \kappa(P)$ de $\varphi$ au-dessus de $P$ est g\'eom\'etriquement r\'eguli\`ere sur le corps $\kappa(P)$ ( c'est-\`a-dire si l'anneau $B\otimes_A \k$ est r\'egulier pour toute extension finie $\k$ de $\kappa(P)$) (cf. \cite{Ma}).
\begin{theoreme}\cite{Ar}\cite{P-P}
Soit $(A,\,\m)$ un couple local n{\oe}th\'erien. Alors si ce couple v\'erifie la propri\'et\'e d'approximation, alors il v\'erifie la propri\'et\'e d'approximation forte.
\end{theoreme}
Ce deuxi\`eme th\'eor\`eme n'est pas vrai dans le cas g\'en\'eral. M. Spivakovsky a donn\'e un exemple de paire hens\'elienne  v\'erifiant la PA et donn\'e un polyn\^ome qui n'admet pas de fonction de Artin \cite{Spi1}.\\
Dans le cas d'un couple local, la fonction de Artin d'un id\'eal $(f)$ de $A[X]$ est une mesure de la non-lissit\'e du morphisme $A\lgw A[X]/(f)$, celle-ci \'etant \'egale \`a l'identit\'e quand ce morphisme est lisse.\\
\\
Le but de cet article est d'utiliser le lemme d'Artin-Rees \cite{Ma} et le th\'eor\`eme d'Izumi \cite{I2} \cite{Re2} pour d\'eterminer une certaine classe de polyn\^omes dont les fonctions de Artin sont born\'ees par des fonctions affines. Nous savons qu'en g\'en\'eral ceci est faux et a pour cons\'equence qu'il n'existe pas de th\'eorie d'\'elimination des quantificateurs dans l'anneau des s\'eries en plusieurs variables muni d'un langage de premier ordre de Presburger \cite{Ro}. N\'eanmoins il existe certains cas pour lesquels ce r\'esultat est vrai.\\
Nous utilisons ici le lemme d'Artin-Rees et le th\'eor\`eme d'Izumi \cite{I2} pour \'etudier la fonction de Artin de certains polyn\^omes.\\
Nous commen\c{c}ons par \'enoncer quelques r\'esultats de r\'eduction. Ensuite, dans la troisi\`eme partie, nous citons le cas des syst\`emes d'\'equations lin\'eaires qui d\'ecoule du lemme d'Artin-Rees (th\'eor\`eme \ref{lin}). Nous montrons dans la quatri\`eme partie que le th\'eor\`eme d'Izumi est \'equivalent \`a une majoration des fonctions de Artin d'une certaine famille de polyn\^omes lin\'eaires (proposition \ref{ICLunif} et th\'eor\`eme \ref{artinreesunif}) et en d\'eduisons une version stable du  lemme d'Artin-Rees (th\'eor\`eme \ref{unif}). Nous donnons ensuite diff\'erentes applications de ces deux r\'esultats:\\
 En cinqui\`eme partie, nous montrons que la fonction de Artin de $X_1X_2-X_3X_4$, vu comme polyn\^ome \`a coefficients dans l'anneau des s\'eries formelles en plusieurs variables, n'est pas born\'ee par une fonction affine.\\
En sixi\`eme partie, nous utilisons simultan\'ement le lemme d'Artin-Rees et le th\'eor\`eme d'Izumi pour montrer que les polyn\^omes qui sont de la forme $f\prod_{k=1}^rX_k^{n_k}+\sum_{j=1}^pf_jZ_j$ ont une fonction de Artin born\'ee par une fonction affine, dans le cas o\`u l'anneau de base quotient\'e par l'id\'eal $(f_1,...,\,f_p)$ est r\'eduit (th\'eor\`eme \ref{thmprinc}).\\ 
Enfin, en derni\`ere partie nous montrons que ceci implique que la fonction de Artin de certains polyn\^omes est born\'ee par une fonction affine (propositions \ref{DS1} et \ref{xnf}) et  nous utilisons ces r\'esultats pour  calculer des cl\^otures int\'egrales approch\'ees d'id\'eaux (exemple \ref{ex}).\\
Je tiens \`a remercier M. Hickel et M. Spivakovsky pour leurs conseils et remarques. Je suis gr\'e au premier de m'avoir fait remarquer que le lemme d'Artin-Rees et le th\'eor\`eme d'Izumi \'etaient des cas particuliers de lin\'earit\'e de fonctions de Artin. Je tiens aussi \`a remercier vivement le referee pour ses remarques et sa patience face \`a une premi\`ere version tr\`es p\'enible.\\
\\
Les anneaux consid\'er\'es seront toujours commutatifs et unitaires.
Nous noterons dans la suite $T=(T_{1},..,T_{N})$, $X=(X_{1},..,X_{n})$ et $f=(f_{1},..,f_{p})$. Sauf indication contraire nous noterons $\m$ l'id\'eal maximal de l'anneau local \'etudi\'e quand il n'y aura aucune confusion possible.

\section{R\'eductions}
Nous allons ici \'enoncer quelques lemmes qui nous permettront de nous ramener \`a \'etudier le cas o\`u l'anneau de base est un anneau complet r\'egulier :

\begin{lemme}\label{lecomp} \cite{P-P}
Soit $(A,\I)$ un couple n{\oe}th\'erien v\'erifiant la PA pour l'id\'eal $(f)$ et tel que l'id\'eal de $\wdh{A}[X]$ engendr\'e par $(f)$ admette une fonction de Artin. Alors $(f)$ admet une fonction de Artin et celle-ci est \'egale \`a celle de l'id\'eal de $\wdh{A}[X]$ engendr\'e par $(f)$.\\
\end{lemme}

\textbf{Preuve :}
Soient $(f)\subset A[X]$, $\wdh{\b}$ sa fonction de Artin  vu comme id\'eal de $\wdh{A}[X]$ et $x\in A$ tel que $f(x)\equiv 0 \quad \mod \, \mfk{I}^{\wdh{\b}(i)+1}$.
 Donc  il existe $x'\in \wdh{A}$ tel que $f(x')=0$ et $x'-x\in \mfk{I}^{i+1}$. Comme $A$ v\'erifie la PA pour $(f)$, il existe $\ovl{x}\in A$ tel que $f(\ovl{x})=0$ et $\ovl{x}-x'\in \I^{i+1}$.\\
En combinant cela on a $\ovl{x}\in A$ tel que $f(\ovl{x})=0$ et $x-\ovl{x}\in \I^{i+1}$.\\
Inversement, soit $\b$  la fonction de Artin  de $(f)$ vu comme id\'eal de $A[X]$. Soit $x\in\wdh{A}$ tel que $f(x)\equiv 0 \ \mod \, \mfk{I}^{\b(i)+1}$. Choisissons $x'\in A$ tel que $x-x'\in\I^{\b(i)+1}$. Nous avons alors $f(x')\equiv 0 \ \mod \, \mfk{I}^{\b(i)+1}$. Donc il existe $\ovl{x}\in A$ tel que $f(\ovl{x})=0$ et $x'-\ovl{x}\in \mfk{I}^{i+1}$. D'o\`u $x-\ovl{x}\in \mfk{I}^{i+1}$.$\quad\Box$\\

\begin{lemme}\label{lequot} \cite{P-P}
Soit $(A,\I)$ un couple n{\oe}th\'erien et $I$ un id\'eal de $A$. Soient $(f)$ un id\'eal de $\frac{A}{I}[X]$, $(F)$ un id\'eal de $A[X]$ \'egal \`a $(f)$ modulo $I$ et $(g_{1},...,g_{q})$ un syst\`eme de g\'en\'erateurs de $I$. Posons $$G_{k}=F_{k}+\sum_{j}Y_{kj}g_{j} \quad k=1,..,m.$$
Alors si $(G)$ admet une fonction de Artin, alors $(f)$ admet une fonction de Artin born\'ee par celle de $(G)$.
\end{lemme}

\textbf{Preuve :}
Soient $(f)$, $(F)$ et $(G)$ comme dans l'\'enonc\'e.  Soit $\b$ la fonction de Artin de $(G)$.\\
Soit $x\in \frac{A}{I}$ tel que $f(x)\equiv 0 \quad \mod \, \mfk{I}^{\b(i)+1}\frac{A}{I} $ avec $i\in \mathbb{N}$. Soit $x'$ un rel\`evement de $x$ dans $A$. Alors $F(x')\equiv 0 \quad \mod \, \mfk{I}^{\b(i)+1}+I$, c'est-\`a-dire qu'il existe des $y_{kj}\in A$ tels que $F(x')+\sum_{j}y_{kj}g_{j}\equiv 0 \quad \mod \, \mfk{I}^{\b(i)+1}$. Il existe alors une solution $(\ovl{x},\ovl{y})$ de ce syst\`eme $G=0$ avec $\ovl{x}\equiv x' \quad \mod \, \mfk{I}^{i+1}$. Modulo $I$ cette solution convient. Et donc $(f)$ admet une fonction de Artin born\'ee par celle de $(G)$.$\quad\Box$\\
Nous \'enon\c{c}ons maintenant un lemme utile pour la suite :
\begin{lemme}\label{lequot2}
Soit $F(X_1,...,\,X_n)\in A[X_1,...,\,X_n]$ o\`u $(A,\I)$ est un couple n{\oe}th\'erien. Soit $I$ un id\'eal de $A$, $\{f_1,...,f_p\}$ et $\{g_1,...,g_q\}$ deux  syst\`emes de g\'en\'erateurs de $I$. Alors les fonctions de Artin de $h_1=F(X_1,...,\,X_n)+\sum_jf_jY_j$ et de $h_2=F(X_1,...,\,X_n)+\sum_lg_lZ_l$ sont \'egales.
\end{lemme}
\textbf{Preuve :} Il nous suffit de montrer le r\'esultat quand $q=p+1$, $g_i=f_i$ pour $1\leq i\leq p$ et $g_q=g_{p+1}\in I$ est quelconque. En effet dans ce cas, par induction nous voyons que la fonction de Artin de $h_1$ (et de la m\^eme mani\`ere celle de $h_2$) est \'egale \`a la fonction de Artin de $F(X_1,...,\,X_n)+\sum_lg_lZ_l+\sum_jf_jY_j$. Donc $h_1$ et $h_2$ ont des fonctions de Artin \'egales.\\
Soit $h_1$ comme dans l'\'enonc\'e et 
$$h_2:=F(X_1,...,\,X_n)+\sum_{j=1}^pf_jY_j+fY_{p+1}$$
o\`u $f\in I$. Nous pouvons \'ecrire $f=\sum_jf_ju_j$ o\`u les $u_j$ sont dans $A$. Notons $\b_i$ la fonction de Artin de $h_i$ ($i=1$ et 2).\\
$\bullet$ Montrons tout d'abord que $\b_2(i)\geq \b_1(i)$ pour tout $i\in\N$. Soient $x_1,...,\,x_n$, $y_1,...,\,y_p\in A$ et $i\in\N$ tels que 
$$h_1(x,\,y)=F(x_1,...,\,x_n)+\sum_j^pf_jy_j\in\I^{\b_2(i)+1}.$$
Nous avons $h_2(x,\,y_1,...,\,y_p,\,0)=h_1(x,\,y_1,...,\,y_p)$, donc par d\'efinition de $\b_2$, il existe $n+p+1$ \'el\'ements $\ovl{x}_1,...,\,\ovl{x}_n$, $\ovl{y}_1,...,\,\ovl{y}_p,\,\ovl{y}_{p+1}$ tels que nous ayons $h_2(\ovl{x},\,\ovl{y}_1,...,\,\ovl{y}_p,\,\ovl{y}_{p+1})=0$, et $\ovl{x}_k-x_k\in\I^{i+1}$, $1\leq k\leq n$, $\ovl{y}_j-y_j\in\I^{i+1}$, $1\leq j\leq p$, $\ovl{y}_{p+1}\in\I^{i+1}$. Notons alors $\ovl{\ovl{y}}_j=\ovl{y}_j+u_j\ovl{y}_{p+1}$, $1\leq j\leq p$. Nous avons alors
$h_1(\ovl{x},\,\ovl{\ovl{y}})=0$
et $\ovl{x}_k-x_k\in\I^{i+1}$, $1\leq k\leq n$, $\ovl{\ovl{y}}_j-y_j\in\I^{i+1}$, $1\leq j\leq p$. Donc $\b_2(i)\geq \b_1(i)$ pour tout $i\in\N$.\\
$\bullet$ Inversement, montrons que $\b_2(i)\leq \b_1(i)$ pour tout $i\in\N$. Soient $x_1,...,\,x_n$, $y_1,...,\,y_p,\,y_{p+1}\in A$ et $i\in\N$ tels que 
$$h_2(x,\,y)=F(x_1,...,\,x_n)+\sum_j^pf_jy_j+fy_{p+1}\in\I^{\b_1(i)+1}.$$
Nous avons 
$$h_1(x,\,y_1+u_1y_{p+1},...,\,y_p+u_py_{p+1})=h_2(x,\,y_1,...,\,y_p,\,y_{p+1}).$$
Donc par d\'efinition de $\b_1$, il existe $\ovl{x}_1,...,\,\ovl{x}_n$, $\ovl{y}_1,...,\,\ovl{y}_p$ tels que nous ayons $h_1(\ovl{x},\,\ovl{y}_1,...,\,\ovl{y}_p)=0$, et  $\ovl{x}_k-x_k\in\I^{i+1}$, $1\leq k\leq n$, $\ovl{y}_j-(y_j+u_jy_{p+1})\in\I^{i+1}$, $1\leq j\leq p$. Notons alors $\ovl{\ovl{y}}_j=\ovl{y}_j-u_jy_{p+1}$, $1\leq j\leq p$, et $\ovl{\ovl{y}}_{p+1}=y_{p+1}$. Nous avons $h_2(\ovl{x},\,\ovl{\ovl{y}})=0$, et $\ovl{x}_k-x_k\in\I^{i+1}$, $1\leq k\leq n$, $\ovl{\ovl{y}}_j-y_j\in\I^{i+1}$, $1\leq j\leq p+1$. Donc $\b_2(i)\leq \b_1(i)$ pour tout $i\in\N$, et donc $\b_1=\b_2$.$\quad\Box$\\
\\
Nous rappelons ensuite le th\'eor\`eme de structure  de I.S. Cohen pour les anneaux complets locaux. \cite{Ma}\\

\begin{definition}
Un anneau de Cohen $R$ est un  corps de caract\'eristique 0 ou un anneau de valuation discr\`ete complet dont le corps r\'esiduel a une caract\'eristique $p> 0$ et dont l'id\'eal maximal est engendr\'e par $p.1$.\\
\end{definition}

\begin{theoreme} \cite{Ma}\label{EGA}
Soit $A$ un anneau local n{\oe}th\'erien complet. Alors il existe un unique anneau de Cohen $R$ tel que $A$ soit isomorphe au quotient d'un anneau de s\'eries formelles $R[[T]]$.
\end{theoreme}

\section{Fonction de Artin d'un syst\`eme lin\'eaire et lemme d'Artin-Rees}
%\subsection{Fonction de Artin d'un syst\`eme lin\'eaire}

Nous avons le r\'esultat suivant qui nous donne la forme de la fonction de Artin d'un syst\`eme d'\'equations lin\'eaires et qui montre au passage que dans le cas lin\'eaire, l'existence de la fonction de Artin n'est absolument pas li\'ee \`a la propri\'et\'e hens\'elienne mais au fait que l'anneau de base est n{\oe}th\'erien.

\begin{theoreme}\label{lin}
Soit $$\left(
f_1^1X_1+\cdots+f_n^1X_n,...,\ f_1^pX_1+\cdots+f_n^pX_n
\right)$$
un id\'eal de polyn\^omes lin\'eaires not\'e $(f)$ de $A[X_1,...,\,X_n]$ o\`u $(A,\I)$ est un couple. Alors  l'id\'eal $(f)$ admet une fonction de Artin born\'ee par la fonction $i\lgm i+i_0$ si et seulement si nous avons la version faible du lemme de Artin-Rees suivante :
$$I\cap\ovl{\I^i}\subset\I^{i-i_0}I \text{ pour } i\geq i_0$$
o\`u $I$ est le sous-$A$-module de $A^p$  engendr\'e par les $(f_j^1,...,\,f_j^p)$ pour $1\leq j\leq n$ et $\ovl{\I^i}$ le sous-$A$-module de $A^p$ \'egal \`a $\oplus_{k=1}^p\I^i$ pour tout entier $i$.\\
En particulier, si $(A,\I)$ est un couple n{\oe}th\'erien, $(f)$ admet une fonction de Artin born\'ee par une fonction lin\'eaire. De plus le plus petit $i_0$ tel que $i\lgm i+i_0$ majore la fonction de Artin de $(f)$ ne d\'epend que du $A$-module $I$.
\end{theoreme}

\textbf{Preuve :}
Avoir $I\cap\ovl{\I^{i+1}}\subset \I^{i+1-i_0}I$ pour $i_0$ une constante positive, cela est \'equivalent \`a ce que pour tout $x_1,...,\ x_n \in A$ tels que 
\begin{equation}\label{sys}\left\{\begin{array}{ccc}
f_1^1x_1+\cdots+f_n^1x_n & \in & \I^{i+1}\\
\vdots  &    &  \\
f_1^px_1+\cdots+f_n^px_n & \in & \I^{i+1}\\
\end{array}\right.\end{equation}
il existe $\e_1,...,\ \e_n\in \I^{i+1-i_0}$ tels que 
$$\left\{\begin{array}{ccc}
f_1^1x_1+\cdots+f_n^1x_n & = & f_1^1\e_1+\cdots+f_n^1\e_n\\
\vdots  & \vdots   & \vdots \\
f_1^px_1+\cdots+f_n^px_n & = & f_1^p\e_1+\cdots+f_n^p\e_n\\
\end{array}\right.$$
En posant $\ovl{x}_k=x_k-\e_k$, cela est \'equivalent \`a ce que pour tout $x_1,...,\ x_n \in A$ qui v\'erifient le syst\`eme (\ref{sys}) pr\'ec\'edent, il existe $\ovl{x}_1$,..., $\ovl{x}_n\in A$ tels que
$$\left\{\begin{array}{ccc}
f_1^1\ovl{x}_1+\cdots+f_n^1\ovl{x}_n & = & 0\\
\vdots  &   &  \\
f_1^p\ovl{x}_1+\cdots+f_n^p\ovl{x}_n & = & 0\\
\end{array}\right.$$
et $\ovl{x}_k-x_k\in \I^{i+1-i_0}$. Cette derni\`ere condition est exactement \'equivalente \`a dire que l'id\'eal $(f)$ admet une fonction de Artin born\'ee par $i\lgm i+i_0$.\\
La derni\`ere assertion d\'ecoule du fait que si $A$ est n{\oe}th\'erien nous avons le lemme d'Artin-Rees (cf. \cite{Ma} par exemple).$\quad\Box$\\

\begin{rque}
T. Wang \cite{Wa} a caract\'eris\'e le plus petit $i_0$ de la proposition pr\'ec\'edente, dans le cas o\`u $A=\k[[T_1,...,\,T_N]]$ et $\I$ est son id\'eal maximal, en terme de bases standards.
\end{rque}

\section{Th\'eor\`eme d'Izumi et version stable du lemme d'Artin-Rees}
\subsection{Th\'eor\`eme d'Izumi et majoration stable de la fonction de Artin d'une famille de polyn\^omes lin\'eaires}
Nous donnons ici l'\'enonc\'e d'un th\'eor\`eme d'Izumi que nous interpr\'etons en terme de lin\'earit\'e de la fonction de Artin d'un certain type de polyn\^ome. Nous donnons tout d'abord une d\'efinition :
\begin{definition}
Soit $(R,\I)$ un couple  n{\oe}th\'erien o\`u $R$ est local et $\I$ un id\'eal $\m$-primaire avec $\m$ l'id\'eal maximal de $R$. Nous noterons $\nu_{R,\,\I}$ la fonction 
\`a valeurs dans $\N\cup\{\infty\}$ d\'efinie par 
$$\forall x\in R\backslash\{0\},\ \nu_{R,\,\I}(x)=n \Longleftrightarrow x\in\I^n \text{ et } x\notin \I^{n+1}$$
$$\text{et }\nu_{R,\,\I}(0)=\infty.$$
On appelle cette fonction l'ordre $\I$-adique sur $R$.\\
Soit $I$ un id\'eal de $R$, nous noterons $\nu_{I,\,\I}$ pour $\nu_{R/I,\,\I}$ quand aucune confusion sur $R$ ne sera possible. Dans le cas o\`u $\I=\m$ est l'id\'eal maximal de $R$, nous noterons $\nu_R:=\nu_{R,\,\I}$ et $\nu_I:=\nu_{R/I,\,\I}$ (la derni\`ere notation n'est \`a pas confondre avec la valuation $I$-adique).

\end{definition}
Une telle d\'efinition est licite d'apr\`es le lemme de Nakayama.\\
Soit $R$ un anneau local n{\oe}th\'erien et $\I$ un id\'eal $\m$-primaire de $R$. Il est clair que  nous avons $\nu_{I,\,\I}(gh)\geq \nu_{I,\,\I}(g)+\nu_{I,\,\I}(h)$  $\forall g,h\in R$. Il y a \'egalit\'e si et seulement si $Gr_{\I}\left(R/I\right)$ est int\`egre. Nous dirons que $\nu_{I,\,\I}$ admet une in\'egalit\'e compl\'ementaire lin\'eaire (ICL) si il existe $a$ et $b$ r\'eels tels que 
$$\nu_{I,\,\I}(gh)\leq a(\nu_{I,\,\I}(g)+\nu_{I,\,\I}(h))+b\quad \forall g,h\in R.$$
Nous dirons dans ce cas que $a$ et $b$ sont des constantes apparaissant dans une ICL pour $(R,\,\I)$. Nous pouvons remarquer que si $a$ et $b$ existent, alors n\'ecessairement $a\geq 1$ et $b\geq 0$.\\
Nous avons alors le
\begin{theoreme}\label{Izumi}\cite{I2}
Soit $R$  un anneau local n{\oe}th\'erien. Alors il existe deux constantes $a$ et $b$ telles que
$$\nu_{R}(gh)\leq a(\nu_{R}(g)+\nu_{R}(h))+b\quad \forall g,h\in R\backslash\{ 0\}$$
si et seulement si $R$ est analytiquement irr\'eductible.
\end{theoreme}
Soit $I$ un id\'eal de $A$, un anneau local n{\oe}th\'erien, engendr\'e par $f_1$,..., $f_p$. Nous notons $R:=A/I$. Notons alors $i_I$ le plus petit entier tel que $i\lgm i+i_I$ majore la fonction de Artin de $f_1X_1+\cdots+f_pX_p\in A[X]$. Pour tout $x\in A$, notons $\b_x$ la fonction de Artin de $xX_0+f_1X_1+\cdots+f_pX_p$. Nous avons alors la

\begin{proposition}\label{ICLunif}
Avec les notations pr\'ec\'edentes, nous avons :
\begin{enumerate}
\item[(i)]
Si $R$ admet une ICL avec les coefficients $a$ et $b$, alors, pour tout $x\in A$, nous avons la majoration uniforme suivante :
$$\forall i\in\N\quad \b_x(i)\leq ai+a\nu_{I}(x)+ai_I+b.$$
\item[(ii)]
Si nous avons une majoration  uniforme de la fonction $\b_x$ par une fonction de la forme $i\lgm ai+c\nu_{I}(x)+b$, avec $a+c\geq1$,  alors le polyn\^ome $XY+\sum_kf_iX_i\in A[X,\,Y,\,X_1,...,\,X_p]$ admet une fonction de Artin born\'ee par la fonction $i\lgm (a+c)(i+i_I)+\max(b,\,i_I)$, et de plus l'id\'eal $I$ est soit premier, soit $\m$-primaire.
\item[(iii)]
Si le polyn\^ome $XY+\sum_kf_iX_i$ admet une fonction de Artin born\'ee par la fonction $i\lgm ai+b$ et si $I$ est premier alors $R$ admet une ICL
$$\nu_{I}(gh)\leq a(\nu_{I}(g)+\nu_{I}(h))+b\quad \forall g,h\in R.$$
\end{enumerate}
\end{proposition}

\textbf{Preuve :}
Montrons (i) :\\
Soient $x_0,\,x_1,...,\,x_p\in A$ tels que 
$$xx_0+f_1x_1+\cdots+f_px_p\in\m^{ai+a\nu_{I}(x)+ai_I+b+1}.$$
Nous avons donc $\nu_{I}(xx_0)\geq ai+a\nu_{I}(x)+ai_I+b+1$. D'o\`u
$$a(\nu_{I}(x)+\nu_{I}(x_0))+b\geq ai+a\nu_{I}(x)+ai_I+b+1$$
$$\nu_{I}(x_0)\geq i+i_I+1.$$
Nous avons donc $x_0=\sum_k f_kz_k+x'_0$ avec $\ord(x'_0)\geq i+i_I+1$, ce qui implique que 
$$\sum_{k=1}^p f_k(x_k+xz_k)\in \m^{i+i_I+1}$$
car $a\geq 1$.
Il existe donc, par d\'efinition de $i_I$, des $t_k\in A$ qui v\'erifient
$$\forall k\geq 1\quad t_k\in x_k+xz_k+\m^{i+1}$$
$$\text{et } \sum_{k=1}^p f_kt_k=0\, .$$
Nous posons alors $\ovl{x}_0=\sum_k f_kz_k$ et $\ovl{x}_k=t_k-xz_k$ pour $k\geq 1$. Nous avons alors 
$$x\ovl{x}_0+f_1\ovl{x}_1+\cdots+f_p\ovl{x}_p=0\text{ et } \forall k\ \ovl{x}_k-x_k\in\m^{i+1}.$$
Donc $\b_x(i)\leq ai+a\nu_{I}(x)+ai_I+b$ pour tout $i$ dans $\N$.\\

Montrons maintenant (ii) :\\
Nous allons tout d'abord montrer la majoration de la fonction de Artin annonc\'ee, puis nous montrerons que $I$ est soit premier, soit $\m$-primaire.\\
Soient $a$, $b$ et $c$ comme dans l'\'enonc\'e. Fixons tout d'abord $i\geq i_I$. Nous allons montrer que la fonction de Artin de $XY+\sum_kf_iX_i$ est major\'ee par la fonction $i\lgm (a+c)i+\max(b,\,i_I)$. Dans ce cas la fonction de Artin du polyn\^ome $XY+\sum_kf_iX_i$ sera major\'ee par $i\lgm (a+c)(i+i_I)+\max(b,\,i_I)$ comme voulue.\\
Soit $i\geq i_I$ et Soient $x$, $y$, $x_1$,...,\, $x_p$ tels que 
\begin{equation}\label{dep}xy+f_1x_1+\cdots+f_px_p\in\m^{(a+c)i+\max(b,\,i_I)+1}.\end{equation}
Nous allons distinguer deux cas, selon que $x$ et $y$ sont tous les deux dans $I+\m^{i+1}$ ou non.
\begin{enumerate}
\item
Supposons que $x$ et $y$ sont dans $I+\m^{i+1}$, c'est-\`a-dire qu'il existe des $z_{1,j}$ et des $z_{2,j}$ tels que $x-\sum_jf_jz_{1,j}\in\m^{i+1}$ et  $y-\sum_jf_jz_{2,j}\in\m^{i+1}$. En multipliant ces deux termes nous voyons que
$$xy-x\sum_jf_jz_{2,j}-y\sum_jf_jz_{1,j}+\sum_jf_jz_{1,j}\sum_jf_jz_{2,j}\in\m^{2i+1}.$$
D'apr\`es cette relation et la relation (\ref{dep}), on est ramen\'e \`a 
$$\sum_jf_j(x_j+yz_{1,j}+xz_{2,j}-\sum_lf_lz_{1,l}z_{2,j})\in\m^{\min(2i,\,(a+c)i+i_I)+1}$$
Par d\'efinition de $i_I$, il existe donc des $t_j$ tels que $\sum_jf_jt_j=0$ et  $$t_j-\left(x_j+xz_{2,j}+yz_{1,j}-\sum_lf_lz_{1,l}z_{2,j}\right)\in\m^{\min(2i,\,(a+c)i+i_I)-i_I+1}\subset\m^{i+1}.$$ Nous posons alors
$$\ovl{x}=\sum_jf_jz_{1,j},\ \ovl{y}=\sum_jf_jz_{2,j}$$
$$\text{et }\ovl{x}_j=t_j-\left(\ovl{x}z_{2,j}+\ovl{y}z_{1,j}-\sum_lf_lz_{1,l}z_{2,j}\right)=-\sum_lf_lz_{2,l}z_{1,j}.$$
Nous avons donc 
$$\ovl{x}\ovl{y}+\sum_jf_j\ovl{x}_j=0,$$
$$\text{et }\ovl{x}-x,\,\ovl{y}-y \text{ et }x_j-\ovl{x}_j\in\m^{i+1} \text{ pour tout }j.$$
\item
Supposons maintenant que $x\in I+\m^{k+1}$ et $x\notin I+\m^{k+2}$ avec $k<i$. Notons 
$$x=\sum_jf_jz_{1,j}+x'$$
avec $\nu_A(x')=k+1$ et $x'\notin I+\m^{\nu_A(x')+1}$. En particulier, nous voyons que $\nu_I(x)=\nu_I(x')=k+1$. Nous avons alors
$$x'y+\sum_jf_j(x_j+yz_{1,j})\in\m^{(a+c)i+\max(b,\,i_I)+1}.$$
Ou encore
$$x'y+\sum_jf_jx'_j\in\m^{(a+c)i+\max(b,\,i_I)+1}$$
avec $x'_j=x_j+yz_{1,j}$.\\
La fonction de Artin de  $x'Y+\sum_k f_kX'_k\in A[Y,X'_1,...,X'_n]$ est major\'ee par $$i\lgm ai+c\nu_{I}(x')+b\leq (a+c)i+b.$$
Donc il existe $\ovl{y}\in y +\m^{i+1}$ et $\ovl{x}'_j \in x'_j +\m^{i+1}$ tels que 
$$x'\ovl{y}+\sum_jf_j\ovl{x}'_j=0\ .$$
Posons alors 
$$\ovl{x}_j=\ovl{x}'_j-\ovl{y}z_{1,j}\text{ et }\ovl{x}=x.$$
Nous avons $$\ovl{x}\ovl{y}+\sum_jf_j\ovl{x}_j=(x'+\sum_jf_jz_{1,j})\ovl{y}+\sum_jf_j(\ovl{x}'_j-\ovl{y}z_{1,j})=0$$
et $\ovl{x}-x\in\m^{i+1}$, $\ovl{y}-y\in\m^{i+1}$ et $\ovl{x}_j-x_j\in\m^{i+1}$ pour tout $j$.
\end{enumerate}
Donc pour $i\geq i_I$ la fonction de Artin de $XY+\sum_kf_iX_i$ est born\'ee par la fonction $i\lgm (a+c)i+\max(b,\,i_I)$.\\

Montrons maintenant que $I$ est premier ou $\m$-primaire. Montrons tout d'abord que $I$ n'a qu'un id\'eal premier minimal associ\'e. Supposons le contraire, c'est-\`a-dire que nous avons $I=I_1\cap I_2$ avec $I\neq I_1$ et $I\neq I_2$, o\`u $I_1$ est un id\'eal $P$-primaire avec $P$ premier, et $P$ n'est pas un id\'eal premier associ\'e \`a $I_2$. Soit $x\in I_1\backslash I_1\cap I_2$. Pour tout entier $l$, il existe $\ovl{x}(l)$ tel que $\nu_{A}(\ovl{x}(l))\geq l$ et $x(l)=x+\ovl{x}(l)\notin P$. En effet, si cela n'\'etait pas possible, nous aurions $x+\m^l\subset P$ pour $l\in \N$. Par cons\'equent, comme $x\in P$ et que $P$ est premier, nous avons $\m\subset P$, et donc $\m=P$, ce qui est impossible par hypoth\`ese sur $P$.\\
Choisissons alors $y\in I_2\backslash I_1\cap I_2$. Il existe un entier $k$ tel que $y\notin I_1+\m^k$ car $y\notin I_1$. Nous avons $xy\in I_1I_2\subset I$, il existe donc des $x_j$ tels que
$$x(l)y=xy+\ovl{x}(l)y=-\sum_j f_jx_j+\ovl{x}(l)y.$$
Donc $x(l)y+\sum_j f_jx_j\in\m^l$. Si $\ovl{y},\,\ovl{x}_1,...,\,\ovl{x}_p$ v\'erifient $x(l)\ovl{y}+\sum_j f_j\ovl{x}_j=0$, alors $x(l)\ovl{y}\in I\subset I_1\subset P$. Donc  $\ovl{y}\in I_1$, car $x(l)\notin P$ et $I_1$ et $P$-primaire. Donc $y-\ovl{y}\notin \m^k$.  D'autre part, pour $l$ assez grand (en fait pour $l>\nu_I(x)$), nous avons $\nu_I(x(l))=\nu_I(x)<+\infty$. La fonction de Artin $\b_x$ n'est donc pas major\'ee uniform\'ement par une fonction de $\nu_I(x)$, ce qui est contradictoire avec l'hypoth\`ese, et donc $I$ n'a qu'un id\'eal premier minimal associ\'e.\\
Supposons maintenant que $I$ n'a qu'un id\'eal minimal associ\'e mais que $I$ n'est ni premier ni $\m$-primaire. C'est-\`a-dire $I$ est $P$-primaire, $I\neq P$ et $P\neq\m$. L'id\'eal $P$ est de la forme $(I:y)$ avec $y\notin I$. Soit $x\in P\backslash P\cap I$. Alors $xy\in I$ par d\'efinition de $y$.\\
Pour tout entier $l$, il existe $\ovl{x}(l)$ tel que $\nu_{A}(\ovl{x}(l))\geq l$ et $x(l)=x+\ovl{x}(l)\notin P$. Si cela n'\'etait pas possible, alors, comme pr\'ec\'edemment, nous aurions $P=\m$ ce qui est contraire \`a l'hypoth\`ese donc impossible.\\ 
Il existe  un entier $k$ tel que $y\notin I+\m^k$ car $y\notin I$. Or $xy\in I$, donc il existe  des $x_j$ tels que
$$x(l)y=xy+\ovl{x}(l)y=-\sum_j f_jx_j+\ovl{x}(l)y.$$
Donc $x(l)y+\sum_j f_jx_j\in\m^l$. Comme pr\'ec\'edemment, si $\ovl{y},\,\ovl{x}_1,...,\,\ovl{x}_p$ v\'erifient $x(l)\ovl{y}+\sum_j f_j\ovl{x}_j=0$, alors $x(l)\ovl{y}\in I\subset P$, donc  $\ovl{y}\in I$ car $x(l)\notin P$ et $I$ est $P$-primaire. Donc $y-\ovl{y}\notin \m^k$. Comme pr\'ec\'edemment, la fonction de Artin $\b_x$ n'est donc pas major\'ee uniform\'ement par une fonction de $\nu_I(x)$ et donc $I$ est premier ou $\m$-primaire.\\

Montrons finalement (iii) :\\
Soient $x$, $y$ et $i$ tels que $a(i+1)+b\geq\nu_{I}(xy)\geq ai+b+1$.
C'est-\`a-dire $xy\in I +\m^{ai+b+1}$. Il existe alors des $z_k$ tel que $xy+\sum_k f_kz_k\in \m^{ai+b+1}$. Il existe donc $\ovl{x}$, $\ovl{y}$ et $\ovl{z}_k$ tels que $\ovl{x}\ovl{y}+\sum_k f_k\ovl{z}_k=0$ et $x-\ovl{x}\in\m^{i+1} $, $y-\ovl{y}\in\m^{i+1}$. Comme $I$ est premier, alors soit $\ovl{y}\in I$, soit $\ovl{x}\in I$. D'o\`u soit $\nu_{I}(x)\geq i+1$, soit $\nu_{I}(y)\geq i+1$. C'est-\`a-dire 
$$\text{soit }a\nu_{I}(x)+b\geq \nu_{I}(xy),$$
$$\text{soit } a\nu_{I}(y)+b\geq \nu_{I}(xy).$$
Nous avons donc $$\nu(xy)\leq a\max(\nu_{I}(x),\,\nu_{I}(y))+b\leq a(\nu_{I}(x)+\nu_{I}(y))+b.$$
D'o\`u le r\'esultat.$\quad\Box$

\begin{rque}
La preuve de ii) pr\'ec\'edente nous montre en fait que, si $I$ n'est ni premier ni $\m$-primaire, nous n'avons aucune majoration uniforme de $\b_x$ par une fonction de $\nu_I(x)$ (m\^eme non affine).\\
\end{rque}

\subsection{Version stable du lemme d'Artin-Rees}
Nous avons en fait le r\'esultat suivant d\^u \`a Rees \cite{Re2} qui est un peu plus fort que celui d'Izumi:
\begin{theoreme}\cite{Re2}
Soit $R$  un anneau local n{\oe}th\'erien. Alors $R$ est analytiquement irr\'eductible si pour au moins un id\'eal $\I$ $\m$-primaire, et seulement si pour tout id\'eal $\I$ $\m$-primaire, il existe deux constantes $a$ et $b$ telles que
$$\nu_{R,\,\I}(gh)\leq \nu_{R,\,\I}(g)+a\nu_{R,\,\I}(h)+b\quad \forall g,h\in R\backslash\{ 0\} \ .$$
\end{theoreme}
Nous en d\'eduisons le
\begin{theoreme}\label{artinreesunif}Soient $A$ un anneau local n{\oe}th\'erien, $I$ un id\'eal de $A$ et $\I$ un id\'eal $\m$-primaire de $A$ o\`u $\m$ est l'id\'eal maximal de $A$, tels que $A/I$ soit analytiquement irr\'eductible. 
 Alors pour tout $x\in A$, nous avons la majoration uniforme suivante :
$$\forall i\in\N\quad \b_x(i)\leq i+a\nu_{I,\,\I}(x)+i_I+b$$
o\`u $\b_x$ est la fonction de Artin de $xX_0+f_1X_1+\cdots+f_pX_p$ pour le couple $(A,\,\I )$.
\end{theoreme}
\textbf{Preuve :}
Soient $x_0,\,x_1,...,\,x_p\in A$ tels que 
$$xx_0+f_1x_1+\cdots+f_px_p\in\I^{i+a\nu_{I,\,\I}(x)+i_I+b+1}.$$
Nous avons donc $\nu_{I,\,\I}(xx_0)\geq i+a\nu_{I,\,\I}(x)+i_I+b+1$. D'o\`u
$$a\nu_{I,\,\I}(x)+\nu_{I,\,\I}(x_0)+b\geq i+a\nu_{I,\,\I}(x)+i_I+b+1$$
$$\nu_{I,\,\I}(x_0)\geq i+i_I+1.$$
Nous avons donc $x_0=\sum_k f_kz_k+x'_0$ avec $\nu_{A,_,\I}(x'_0)\geq i+i_I+1$, ce qui implique que 
$$\sum_{k=1}^p f_k(x_k+xz_k)\in \I^{i+i_I+1}.$$
Il existe donc, par d\'efinition de $i_I$, des $t_k\in A$ qui v\'erifient
$$\forall k\geq 1\quad t_k\in x_k+xz_k+\I^{i+1}\,\text{ et } \sum_{k=1}^p f_kt_k=0\, .$$
Nous posons alors $\ovl{x}_0=\sum_k f_kz_k$ et $\ovl{x}_k=t_k-xz_k$ pour $k\geq 1$. Nous avons alors 
$$x\ovl{x}_0+f_1\ovl{x}_1+\cdots+f_p\ovl{x}_p=0\text{ et } \forall k\ \ovl{x}_k-x_k\in\I^{i+1}\, .\quad\Box$$\\
Nous pouvons alors formuler  une version stable du lemme d'Artin-Rees :

\begin{theoreme}\label{unif}
Soient $A$ un anneau n{\oe}th\'erien, $\I$ un id\'eal $P$-primaire de $A$ avec $P$ premier et $I\subset P$ un id\'eal de $A$ tel que $A_P/IA_P$ soit analytiquement irr\'eductible. Supposons que 
\begin{enumerate}
\item[i)] $\forall k\geq 1,\ \I^kA_P\cap A=\I^k,$
\item[ii)] $\forall k\geq 1,\, \forall x\in P,\ ((x)+I)\I^kA_P\cap A=((x)+I)\I^k.$
\end{enumerate}
Alors il existe $a\geq 1$ et $b\geq 0$ tels que nous ayons la version faible d'Artin-Rees uniforme suivante
$$\left((x)+I\right)\cap \I^{i+a\nu_{I,\,\I}(x)+b}\subset \left((x)+I\right)\I^i\quad  \forall x\in P\ \forall i\in\N.$$
\end{theoreme}
\textbf{Preuve :} D'apr\`es i), les ordres $\nu_{A,\,\I}$ et $\nu_{A_P,\I A_P}$ sont \'egaux. D'apr\`es le th\'eor\`eme pr\'ec\'edent et le th\'eor\`eme \ref{lin}, il existe $a$ et $b$ tels que 
$$\left((x)+I\right)A_P\cap \I^{i+a\nu_{I,\,\I}(x)+b}A_P\subset \left((x)+I\right)\I^iA_P\quad  \forall x\in PA_P\ \forall i\in\N$$
car $A_P/IA_P$ est  analytiquement irr\'eductible.
Choisissons $x\in P$ et $i\in\N$, nous avons alors
$$\left((x)+I\right)\cap \I^{i+a\nu_{I,\,\I}(x)+b}\subset\left((x)+I\right)A_P\cap \I^{i+a\nu_{I,\,\I}(x)+b}A_P \subset \left((x)+I\right)\I^iA_P\, .$$
Le r\'esultat d\'ecoule alors de l'hypoth\`ese ii).$\quad\Box$

\begin{rque}
Ceci est vrai en particulier si $A$ est local, $P=\m$ est son id\'eal maximal, $\I$ est $\m$-primaire  et $A/I$ est analytiquement irr\'eductible.
\end{rque}
\begin{rque}
Il existe deux versions de ce que l'on appelle lemme d'Artin-Rees uniforme \cite{Hu} et \cite{B-M} qui sont \`a ne pas confondre avec cette version stable.
\end{rque}

\subsection{Exemples}
Nous donnons ici quelques exemples explicites, toujours dans le cas o\`u l'id\'eal $\I$ est l'id\'eal maximal de l'anneau $A$. Nous noterons alors  $\m$ cet id\'eal. Dans la suite, l'anneau $\O_N$ d\'esignera indiff\'eremment l'anneau des s\'eries formelles en $N$ variables sur un corps $\k$ et l'anneau des s\'eries convergentes en $N$ variables sur $\k$ (quand cela a un sens). Nous noterons $\ord$  l'ordre $\m$-adique sur $\O_N$.
\subsubsection{Premier exemple}\label{gradint}
Si l'anneau gradu\'e $Gr_{\m}\frac{A}{I}$ est int\`egre alors $\nu_{A,\,I}$ est une valuation, i.e.
$$\nu_{A,\,I}(gh)= (\nu_{A,\,I}(g)+\nu_{A,\,I}(h))\quad \forall g,h\in A$$
En particulier d'apr\`es le th\'eor\`eme \ref{ICLunif}, la fonction de Artin du polyn\^ome $xX_0+f_1X_1+\cdots+f_pX_p$ (o\`u $I=(f_1,...,\,f_p)$) est born\'ee par une fonction de la forme $i\lgm i+\nu_{A,\,I}(x)+p$.\\
C'est le cas par exemple si $I=(f)$ et $f$ est irr\'eductible et homog\`ene de degr\'e $p$ dans $\O_N$.\\

\subsubsection{Deuxi\`eme exemple}
Nous allons donner tout d'abord le
\begin{lemme}\label{lininiirr}
Soit $L(X_1,...,\,X_n)=f_1X_1+\cdots+f_nX_n \in\O_N[X_1,...,\,X_n]$ avec $\ord(f_1)\leq \ord(f_2)\leq ...\leq \ord(f_n)$. Supposons que les termes de plus bas ordre (termes initiaux) des $f_k$  forment une suite r\'eguli\`ere. Alors $L$ admet une fonction de Artin qui est major\'ee,  pour tout $i\geq 0$, par la  fonction affine $i\lgm i+\ord(f_n)$.
\end{lemme}

\textbf{Preuve :}
Les termes initiaux des $f_k$ formant une suite r\'eguli\`ere, les $f_k$ forment une suite r\'eguli\`ere et nous savons donc que les z\'eros de $L$ sont de la forme 
$$\left(\sum_{k=1}^nf_kz(k,1),...,\,\sum_{k=1}^nf_kz(k,n)\right)$$
avec $z(k,j)=-z(j,k)$ pour tous $k$ et $j$.  En particulier $z(k,k)=0$ pour tout $ k$.\\
Dans la suite, pour tout \'el\'ement $x$ de $\O_N$, nous noterons $x(p)$ le terme homog\`ene de degr\'e $p$ de $x$.\\
Soient $x_1$,..., $x_n\in\O_N$ tels que $f_1x_1+\cdots+f_nx_n \in\m^{i+\ord(f_n)+1}$. Si nous avons $\min_j(\ord(f_jx_j))\geq i+\ord(f_n)+1$, nous posons $\ovl{x}_j=0$ pour tout $j$. Nous avons $L(\ovl{x})=0$ et $x_j-\ovl{x}_j\in\m^{i+1}$ pour tout $j$.\\
Dans le cas contraire, nous allons construire, par r\'ecurrence sur $\min_j(\ord(f_jx_j))$, des \'el\'ements $\ovl{x}_j$, pour tout $j$, tels que $\sum_jf_j\ovl{x}_j=0$ et $\ovl{x}_j-x_j\in\m^{i+\ord(f_n)-\ord(f_j)+1}$ pour tout $j$.\\
Comme $\min_j(\ord(f_jx_j))<i+\ord(f_n)+1$, nous avons 
$$in\left(\sum_{j=1}^n f_{j}(\ord(f_j))x_{j}( \ord(x_j))\right)=0$$
o\`u $in(x)$ d\'esigne le terme initial de $x$ pour $\ord$.
C'est-\`a-dire
$$\sum_{j\in I_1} f_{j}( \ord(f_j))x_{j}( \ord(x_j))=0$$
o\`u $I_1$ est l'ensemble
$$I_1:=\left\{j\in\{1,...,\,n\}\,/\, \ord(f_jx_j)\leq \ord(f_kx_k),\ \forall k\in\{1,...,\,n\}\right\}.$$
Il existe donc des polyn\^omes homog\`enes $z^1(k,j)\in\O_N$ tels que 
$$z^1(k,j)=0 \text{ si } j\notin I_1,\,z^1(k,j)=-z^1(j,k)$$
$$\text{et }x_{j}( \ord(x_j))=\sum_{k=1}^n f_{k}( \ord(f_k))z^1(k,j) \text{ pour tout }j\in I_1$$
car les termes initiaux des $f_j$, o\`u $j\in I_1$, forment une suite r\'eguli\`ere.
Nous posons alors 
$$x^1_j=x_j-\sum_{k=1}^n f_kz^1(k,j)\ \forall j.$$
Nous avons donc $f_1x^1_1+\cdots+f_nx^1_n \in\m^{i+\ord(f_n)+1}$ et $\ord(x^1_j)>\ord(x_j)$ si $j\in I_1$ et $\ord(x^1_j)=\ord(x_j)$ sinon. Nous avons aussi que $$\min_j(\ord(f_jx_j))<\min_j(\ord(f_jx^1_j)).$$
Nous pouvons alors continuer ce processus jusqu'au rang $l$ de mani\`ere \`a  avoir contruit des $x^l_j$ tels que $f_jx^l_j\in\m^{i+\ord(f_n)+1}$ pour tout $j$ avec
$$x^l_j=x_j-\sum_{k=1}^n f_k\ovl{z}(k,j)\text{ tels que } \ovl{z}(k,j)=-\ovl{z}(j,k) \ \forall k,\,j.$$
C'est-\`a-dire qu'il existe $\ovl{x}_j=\sum_{j=1}^n f_k\ovl{z}(k,j)$ tels que 
$$f_1\ovl{x}_1+\cdots+f_n\ovl{x}_n =0$$
$$\text{et } \, \forall j,\ x_j-\ovl{x}_j\in\m^{i+\ord(f_n)-\ord(f_j)+1}\subset \m^{i+1}.
\quad\Box$$
Nous en d\'eduisons le  
\begin{corollaire}\label{coefinit}
Soit $I=(f_1,...,\,f_n)$ un id\'eal de $\O_N$. Si l'id\'eal engendr\'e par les termes initiaux des \'el\'ements de $I$ est premier et d'intersection compl\`ete  alors nous avons l'in\'egalit\'e
$$\nu_{I}(gh)\leq 2(\nu_{I}(g)+ \nu_{I}(h))+3\max_k\{\ord(f_k)\}\quad \forall f,g\in\O_N.$$
\end{corollaire}
\textbf{Preuve :} Soit $f_1,...,\,f_n$ une famille d'\'el\'ements de $I$ dont les termes initiaux forment une suite r\'eguli\`ere et engendrent l'id\'eal des termes initiaux  de $I$. Alors cette famille engendre $I$ en tant qu'id\'eal. Soit $f\in\O_N$ et $f'$ son reste apr\`es division par $I$ (th\'eor\`eme de division d'Hironaka cf. \cite{Hi}). Si $f'=0$, alors $f\in I$ et la fonction de Artin de $fX_0+f_1X_1+\cdots+f_nX_n$ est born\'ee par $i\lgm i+i_I$.\\
Si $f'\neq 0$, alors $\nu_{I}(f)=\nu_{I}(f')$, et la suite form\'ee des termes initiaux des $f_l$ et du terme initial de $f'$ est r\'eguli\`ere. En effet, en notant $in(g)$ le terme initial de $g\in\O_N$, supposons qu'il existe $x\in in(I+(f))$ tel que nous ayons $x\, in(f')=0$ dans $in(I+(f))/(in(f_1,...,f_n))$. Comme $in(I)$ est premier, n\'ecessairement $x\in in(I)$ et donc la suite $(in(f_1),...,\,in(f_n),\,in(f'))$ est r\'eguli\`ere.\\
D'apr\`es le th\'eor\`eme \ref{lin}, la fonction de Artin de $fX_0+f_1X_1+\cdots+f_nX_n$ est \'egale \`a celle de $f'X_0+f_1X_1+\cdots+f_nX_n$, qui est born\'ee, d'apr\`es le lemme \ref{lininiirr}, par 
$$i\lgm i +\max\{\ord(f'),\,\ord(f_k)\}\leq i+\ord(f')+\max_k\{\ord(f_k)\}.$$ En utilisant alors le (ii) de la proposition \ref{ICLunif} ($a=c=1$ et $b=\max_k\{\ord(f_k)\}$), nous voyons que le polyn\^ome $XY+\sum_kf_kZ_k$ admet une fonction de Artin born\'ee par la fonction
$i\lgm 2(i+i_I)+\max\left(\max_k\{\ord(f_k)\},\,i_I\right)$ ($i_I$ est une constante telle que $i\lgm i+i_I$ majore la fonction de Artin de $\sum_kf_kZ_k$). Comme $i_I\leq \max_k\{\ord(f_k)\}$ d'apr\`es le lemme \ref{lininiirr}, le polyn\^ome $XY+\sum_kf_kZ_k$ admet une fonction de Artin born\'ee par la fonction $i\lgm 2i+3\max_k\{\ord(f_k)\}$. En utilisant alors le (iii) de la proposition \ref{ICLunif}, nous d\'eduisons le r\'esultat.$\quad\Box$

\subsubsection{Troisi\`eme exemple}
Soit $f=T_1^2+g(T_2,\,T_3)\in\O_3$ avec $g(0,\,0)=0$. Alors d'apr\`es \cite{I2}, $(f)$ admet une ICL avec les coefficients $1$ et $\ord(g)-2$ si $\ord(g)$ est impair. En utilisant le (i) de la proposition \ref{ICLunif} ($a=1$ et $b=\ord(g)-2$), nous voyons que  la fonction de Artin de $xX_0+fX_1$ est born\'ee par 
$$i\lgm i+\nu_{(f),\m}(x)+\ord(g)-2+i_{(f)}$$
o\`u $i_{(f)}$ est tel que $i\lgm i+i_{(f)}$ majore la fonction de Artin de $fX_1$. En particulier nous pouvons choisir $i_{(f)}=\ord(f)=2$. Donc la fonction de Artin de $xX_0+fX_1$ est born\'ee par 
$$i\lgm i+\nu_{(f),\m}(x)+\ord(g).$$

\subsubsection{Quatri\`eme exemple}
Nous allons donner une ICL dans le cas o\`u $f=T_1^k+g\in\O_N$  avec $\ord(g)=k+1$ et $T_1$ ne divisant pas le terme initial de $g$. Nous avons tout d'abord le

\begin{lemme}
Soit $f=T_1^k+g$ avec $\ord(g)=k+1$  et $T_1$ ne divisant pas le terme initial de $g$. Alors pour tout $h$ la fonction de Artin de $fX+hY$ est born\'ee par 
$$i\lgm i+\max\{k,\,\nu_{f,\m}(h)+1\}.$$
\end{lemme}

\textbf{Preuve :}
Soit $h=af+h_0T_1^l+\sum_{j\geq 1}h_j$ avec $l<k$ et $T_1$ ne divisant pas $h_0$, et les $h_j$ sont homog\`enes de degr\'e $j>\ord(h_0)+l$ et ne sont pas divisibles par $T_1^k$. Notons $h'=h_0T_1^l+\sum_{j\geq 1}h_j$.\\
Soient $x$ et $y$ tels que $fx+hy\in\m^{i+\max\{k,\,\nu_{f,\m}(h)\}+2}$. Nous avons donc 
$$f(x+ay)+h'y\in\m^{i+\max\{k,\,\nu_{f,\m}(h)\}+2}.$$
Nous pouvons faire le changement de variables $X=X+aY$, $Y=Y$ et supposer que $h=h'$.\\
Notons $x_{j}$ le terme homog\`ene de degr\'e $j$ dans l'\'ecriture de $x$ (idem pour $y$). Si $\ord(x)\geq i+\max\{k,\,\nu_{f,\m}(h)+1\}-k+1\geq i+1$, nous posons $\ovl{x}=\ovl{y}=0$. Nous avons bien  $\ovl{x}-x$, $\ovl{y}-y\in\m^{i+1}$ et $f\ovl{x}+h'\ovl{y}=0$.\\
Autrement nous avons
$$T_1^kx_{\ord(x)}+h_0T_1^ly_{\ord(y)}=0$$
$$T_1^kx_{\ord(x)+1}+in(g)x_{\ord(x)}+h_0T_1^ly_{\ord(y)+1}+h_1y_{\ord(y)}=0.$$\\
La premi\`ere \'equation nous donne que $T_1^{k-l}$ divise $y_{\ord(y)}$. La seconde \'equation nous donne alors que $T_1^{min\{l,\,k-l\}}$ divise $x_{\ord(x)}$.\\
Si $l\leq k-l$ alors nous avons $x_{\ord(x)}=h_0T_1^lz_0$ et  $y_{\ord(y)}=-T_1^kz_0$. Nous posons alors $x(1)=x-hz_0$ et $y(1)=y+fz_0$. Nous avons $\ord(x(1))>\ord(x)$ et $\ord(y(1))>\ord(y)$.\\
Si $l>k-l$, la premi\`ere \'equation nous donne que $T_1^{2(k-l)}$ divise $y_{\ord(y)}$ et la seconde que $T_1^{min\{l,\,2(k-l)\}}$ divise $x_{\ord(x)}$.\\
Par induction nous pouvons continuer cette proc\'edure jusqu'au rang $p$ tel que $l\leq p(k-l)$ et tel que $T_1^{min\{l,\,p(k-l)\}}=T_1^l$ divise $x_{\ord(x)}$. Il existe donc $z_0$ tel que $x_{\ord(x)}=h_0T_1^lz_0$ et  $y_{\ord(y)}=-T_1^kz_0$. Nous posons alors $x(1)=x-hz_0$ et $y(1)=y+fz_0$. Nous avons $\ord(x(1))>\ord(x)$ et $\ord(y(1))>\ord(y)$.\\
Nous recommen\c{c}ons alors la proc\'edure pr\'ec\'edente et nous construisons ainsi $z$ tel que $\ord(x-hz)\geq i+\max\{k,\,\nu_{f,\m}(h)\}-k+1\geq i+1$. Nous posons alors $\ovl{x}=hz$ et $\ovl{y}=-fz$. Clairement $\ovl{x}-x$, $\ovl{y}-y\in\m^{i+1}$ et $f\ovl{x}+h'\ovl{y}=0$. $\quad\Box$\\
\\
D'apr\`es le (ii) la proposition \ref{ICLunif} (avec $a=c=1$ et $b=k$), nous voyons donc que le germe d'hypersurface d\'efini par $f=T_1^k+g=0$ avec $\ord(g)=k+1$ et $\pgcd(T_1,\, in(g))=1$ admet une ICL :
$$\nu_{(f),\m}(gh)\leq 2(\nu_{(f),\m}(g)+\nu_{(f),\m}(h))+3k\quad \forall g,h\in\O_N.$$

\section{Etude de la fonction de Artin de $X_1X_2-X_3X_4$}
Nous donnons ici un exemple de polyn\^ome dont la fonction de Artin n'est pas born\'ee par une fonction affine. L'id\'ee est d'utiliser le fait que la fonction de Artin de $X_1X_2-(T_1T_2-T_3^i)X_4\in\O_N[X_1,\,X_2,\,X_4]$ pour $N\geq3$ est la fonction $k\lgm ik-1$ (cf. exemple 5.6 (iv) de \cite{I2}) et que tout \'el\'ement \'egal \`a $T_1T_2-T_3^i$ modulo $\m^{i+1}$ est toujours irr\'eductible. Ce polyn\^ome peut \^etre alors vu comme une ``sp\'ecialisation'' du polyn\^ome $X_1X_2-X_3X_4\in\O_N[X_1,\,X_2,\,X_3,\,X_4]$.

\begin{theoreme}La fonction de Artin du polyn\^ome 
$$X_1X_2-X_3X_4\in\O_N[X_1,\,X_2,\,X_3,\,X_4]$$
 est  born\'ee inf\'erieurement par la fonction $i\lgm i^2-1$ si $N\geq 3$.
\end{theoreme}
Nous savions d\'ej\`a qu'en g\'en\'eral une fonction de Artin n'\'etait pas born\'ee par une fonction affine (cf. \cite{Ro}). L'exemple \'etudi\'e ici correspond \`a une singularit\'e isol\'ee d'hypersurface, dont la fonction de Artin-Greenberg a d\'ej\`a \'et\'e \'etudi\'ee (cf. \cite{L-J}).\\
\\
\textbf{Preuve :}
Appelons $P$ le polyn\^ome $X_1X_2-X_3X_4$ et fixons un entier $i\in\N$ quelconque. 
Notons $x_1(i):=T_1^i$, $x_2(i):=T_2^i$ et $x_3(i):=T_1T_2-T_3^i$. Nous avons
$$x_1(i)x_2(i)=\left(x_3(i)+T_3^i\right)^i=x_3(i)x_4(i)+T_3^{i^2}$$
avec $x_4(i)$ bien choisi.
Nous avons donc 
$$P(x_1(i),\,x_2(i),\,x_3(i),\,x_4(i))\in\m^{i^2}.$$
Supposons que nous ayons $x_1,\,x_2,\,x_3$ et $x_4$ tels que
$P(x_1,\,x_2,\,x_3,\,x_4)=0$, alors deux cas peuvent se produire :\\
(1) soit $x_3-x_3(i)\in\m^{i+1}$.  Alors $x_3$ est irr\'eductible. En effet, supposons le contraire, c'est-\`a-dire qu'il existe $x$ et $y$ tels que $xy=x_3$. Alors $xy-x_3(i)\in\m^{i+1}$, ce qui est impossible. En effet, d'apr\`es le lemme \ref{irr} dont nous donnons la preuve \`a la fin, la fonction de Artin du polyn\^ome $XY-x_3(i)$ vaut $i$, et cela impliquerait que $x_3(i)$ est r\'eductible, ce qui est clairement faux. Donc soit $x_1\in (x_3)$, soit $x_2\in (x_3)$ car $(x_3)$ est irr\'eductible et $\O_N$ est factoriel. Or 
$$\sup_{f\in\O_N}\big(\ord(x_1(i)-fx_3)\big)=\sup_{f\in\O_N}\big(\ord(x_2(i)-fx_3)\big)=i$$
car $x_1(i)-fx_3=x_1(i)-fx_3(i)$ modulo $\m^i$ et ce dernier terme est non nul modulo $\m^i$, le terme initial de $x_1(i)$ n'\'etant pas divisible par $T_1T_2$ (idem pour $x_2(i)$).\\
(2) soit $\ord\big(x_3-x_3(i)\big)\leq i$.\\
Dans tous les cas nous avons
$$\sup\left(\min_{j=1,..,4}\left(\ord(x_j(i)-x_j)\right)\right)\leq i$$
o\`u la borne sup\'erieure est prise sur tous les 4-uplets $(x_1,\,x_2,\,x_3,\,x_4)$ tels que $P(x_1,\,x_2,\,x_3,\,x_4)=0$. La fonction de Artin de $P$ est donc minor\'ee par la fonction $i\lgw i^2-1$.$\quad\Box$\\
Nous donnons maintenant la preuve du lemme utilis\'e :
\begin{lemme}\label{irr}
La fonction de Artin du polyn\^ome $XY-x_3(i)\in\O_N[X,\,Y]$ est la fonction constante \'egale \`a $i$.
\end{lemme}
\textbf{Preuve :}
Soient $x$ et $y$ dans $\O_N$, non inversibles, tels que $xy-x_3(i)\in\m^{i+1}$. Ecrivons 
$$x=\sum_{j=1}^{i+1}x_j\ \text{ et }\ y=\sum_{j=1}^{i+1}y_j$$
o\`u $x_j$ (resp. $y_j$) est le terme homog\`ene d'ordre $j$ dans l'\'ecriture de $x$ (resp. de $y$). Quitte \`a intervertir $x$ et $y$, nous avons n\'ecessairement $x_1=aT_1$ et $y_1=a^{-1}T_2$. Nous allons montrer par induction, que pour tout $j\in\{1,...,\,i-2\}$, $x_j\in (T_1)$ et $y_j\in (T_2)$. Supposons que ceci soit vrai pour $j\in\{1,...,\,n-1\}$ avec $n<i-1$. Le terme homog\`ene d'ordre $n+1$ de $xy$ est nul car $n+1<i$. Nous avons alors 
$$aT_1y_{n}+a^{-1}T_2x_n+\sum_{j=2}^{n-1}x_jy_{n+1-j}=0$$
Par hypoth\`ese de r\'ecurrence, $\sum_{j=2}^{n-1}x_jy_{n+1-j}\in (T_1T_2)$. Par factorialit\'e de $\O_N$, nous voyons donc que $y_n\in (T_2)$ et $x_n\in (T_1)$.\\
Le terme homog\`ene d'ordre $i$ de $xy$ est donc \'egal \`a
$$aT_1y_{i-1}+a^{-1}T_2x_{i-1}+\sum_{j=2}^{i-2}x_jy_{i-j}.$$
Or ce terme appartient \`a l'id\'eal engendr\'e par $T_1$ et $T_2$. Il ne peut donc pas \^etre \'egal \`a $T_3^i$. Il n'existe donc pas de tels $x$ et $y$, d'o\`u le r\'esultat.$\quad\Box$
\section{Fonction de Artin d'un mon\^ome}
Nous allons utiliser  ici les r\'esultats pr\'ec\'edents  pour montrer  que la fonction de Artin de certains polyn\^omes, en particulier des mon\^omes, est born\'ee par une fonction affine, dans le cas o\`u l'anneau de base est r\'eduit ou analytiquement irr\'eductible.
Nous avons tout d'abord le r\'esultat suivant qui est un corollaire direct de la proposition \ref{ICLunif} :
\begin{corollaire}\label{Izlin}
Soit $$g(X,\,Y,\,Z_j):=XY+\sum_{j=1}^pf_jZ_j$$
avec $I=(f_1,...,\,f_p)$ un id\'eal propre de $A$ n{\oe}th\'erien tel que $A/I$ soit analytiquement irr\'eductible. Alors $g$ admet une fonction de Artin major\'ee par une fonction affine.
\end{corollaire}

Nous donnons ensuite une g\'en\'eralisation du corollaire \ref{Izlin} :

\begin{theoreme}\label{thmprinc}
Soient $A$ un anneau local  n{\oe}th\'erien et $I=(f_j)$ un id\'eal de $A$ tels que  $A/I$ soit analytiquement irr\'eductible ou tels que  $A/I$ soit r\'eduit et $A$ v\'erifie la PA. Alors tout polyn\^ome  \`a coefficients dans $A$ de la forme $f\prod_{k=1}^rX_k^{n_k}+\sum_{j=1}^pf_jZ_j$ admet une fonction de Artin major\'ee par une fonction lin\'eaire. 
\end{theoreme}

\begin{rque}
Le th\'eor\`eme pr\'ec\'edent est vrai en particulier pour un mon\^ome vu comme polyn\^ome \`a coefficients dans un anneau analytiquement irr\'eductible ou r\'eduit et v\'erifiant la PA.
\end{rque}

\textbf{Preuve :} Notons $g(X_k,\,Z_j)=f\prod_{k=1}^rX_k^{n_k}+\sum_{j=1}^pf_jZ_j$.\\
\\
\textbf{Premi\`ere \'etape :} Nous allons d'abord nous ramener au cas o\`u $I=(0)$, c'est-\`a-dire au cas o\`u $g$ est un mon\^ome. Nous notons $\ovl{g}(X_k)$ le polyn\^ome $f\prod_{k=1}^rX_k^{n_k}\in A/I[X_k]$ et supposons que ce polyn\^ome admette une fonction de Artin born\'ee par une fonction affine $a\lgm ai+b$. Soient $x_1,...,\,x_r,\,z_1,...,\,z_p$ tels que $g(x_k,\,z_j)\in\m^{i+1}$. Alors $\ovl{g}(x_k)\in\m^{i+1}$ et donc il existe $\ovl{x}_k\in A$ tel que $\ovl{g}(\ovl{x_k})=0$ dans $A/I$ et $\ovl{x}_k-x_k\in\m^{\frac{i-b}{a}}$. Donc il existe des $z'_j$ tels que $f\prod_{k=1}^r\ovl{x}_k^{n_k}=\sum_jf_jz'_j$ dans $A$. D'o\`u $\sum_jf_j(z_j+z'_j)\in\m^{\frac{i-b}{a}}$ et d'apr\`es Artin-Rees (th\'eor\`eme \ref{lin}) il existe des $t_j$ tels que $\sum_jf_jt_j=0$ et $t_j-(z_j+z'_j)\in\m^{\frac{i-b}{a}-i_0}$ o\`u $i_0$ ne d\'epend que de $I$. Nous posons alors $\ovl{z}_j=t_j-z'_j$ pour tout $j$. Nous avons alors $g(\ovl{x}_k,\,\ovl{z}_j)=0$, et $\ovl{x}_k-x_k\in\m^{\frac{i-b}{a}}$ et $\ovl{z}_j-z_j\in\m^{\frac{i-b}{a}-i_0}$ pour tous $k$ et $j$. Il nous suffit donc de montrer que $\ovl{g}$ admet une fonction de Artin born\'ee par une fonction affine.\\
\\
\textbf{Deuxi\`eme \'etape :}
Nous allons nous ramener au cas o\`u $f=1$. Nous avons $f\prod_{k=1}^rx_k^{n_k}=0$ si et seulement si $\prod_{k=1}^rx_k^{n_k}\in ((0):f)$. De plus si nous avons $f\prod_{k=1}^rx_k^{n_k}\in\m^{i+1}$, alors d'apr\`es Artin-Rees, il existe $i_0$ qui ne d\'epend que de $((0):f)$, tel que $\prod_{k=1}^rx_k^{n_k}\in((0):f)\m^{i-i_0+1}$. Donc montrer que le polyn\^ome $f\prod_{k=1}^rX_k^{n_k}\in A[X_k]$ admet une fonction de Artin born\'ee par une fonction affine revient \`a montrer que $\prod_{k=1}^rX_k^{n_k}\in A/((0):f)[X_k]$ admet une fonction de Artin born\'ee par une fonction affine.\\
Nous pouvons remarquer que si $A$ est r\'eduit et si $x^k\in((0):(f))$ alors $fx^k=0$ et donc $xf=0$ et $x\in((0):f)$, d'o\`u $((0):f)$ est radical et $A/((0):f)$ est r\'eduit.\\
De m\^eme nous pouvons remarquer que si $A$ est analytiquement irr\'eductible alors $A$ est int\`egre et donc $((0):f)=(0)$. Donc  $A/((0):f)=A$ est analytiquement irr\'eductible.\\
\\
\textbf{Troisi\`eme \'etape :} Nous allons traiter le cas o\`u $A/I$ est analytiquement irr\'eductible. Supposons que $f=1$ et $I=(0)$. Soit $i\in\N$ et soient $x_1,...,\,x_r$ tels que $g(x_k)\in\m^{i+1}$. Alors nous avons
$$\nu_I(\prod_{k=1}^rx_k^{n_k})\geq i+1$$
$$\text{et } a\left(\nu_I(\prod_{k=1}^{r-1}x_k^{n_k})+\nu_I(x_r^{n_r})\right)+b\geq i+1$$
o\`u $a$ et $b$ sont les constantes d'une ICL v\'erifi\'ee par $I$. Par r\'ecurrence sur $r$ il existe $k_0\in\{1,...,r\}$ tel que 
$$\nu_I(x_{k_0}^{n_{k_0}})\geq \left\lfloor\frac{i-b'}{a'}\right\rfloor+1$$
pour $a'$ et $b'$ des constantes ind\'ependantes des $x_k$ et de $i$ et o\`u $\left\lfloor c\right\rfloor$ est la partie enti\`ere de $c$. Ensuite si $\nu_I(x^n)\geq \left\lfloor\frac{i-b'}{a'}\right\rfloor+1$, alors par r\'ecurrence sur $n$ nous avons $$\nu_I(x)\geq \left\lfloor\frac{i-b''}{a''}\right\rfloor+1$$
pour $a''$ et $b''$ des constantes ind\'ependantes de $x$ et de $i$. Il suffit alors de poser $\ovl{x}_{k_0}=0$ et $\ovl{x}_{k}=x_k$ pour $k\neq k_0$. Nous avons alors $g(\ovl{x}_k)=0$ et $\ovl{x}_k-x_k\in\m^{\left\lfloor\frac{i-b''}{a''}\right\rfloor+1}$.
Donc le th\'eor\`eme est prouv\'e pour $A/I$ analytiquement irr\'eductible.\\
\\
\textbf{Quatri\`eme \'etape :} Nous allons montrer qu'il suffit, dans le cas o\`u $A$ est  r\'eduit et v\'erifie la PA, de montrer le r\'esultat pour $A$ complet n{\oe}th\'erien et r\'egulier et $I$ radical. Cela d\'ecoule  des lemmes \ref{lecomp} et \ref{lequot}, et du lemme suivant :
\begin{lemme}(\cite{KMPPR}, section 4)
 Soit $A$ un anneau local r\'eduit n{\oe}th\'erien v\'erifiant la PA. Alors $\wdh{A}$ (le compl\'et\'e de $A$ pour la topologie $\m$-adique) est r\'eduit.
\end{lemme}
\textbf{Derni\`ere \'etape :} 
Supposons maintenant que $A$ est complet, n{\oe}th\'erien et r\'egulier et $I$ radical et soient  $i\in\N$ et  $x_1,...,\,x_r,\,z_1,...,\,z_p$ fix\'es tels que $g(x_k,\,z_j)\in\m^{i+1}$. Soit $$I=P_1\cap\cdots\cap P_q$$
la d\'ecomposition primaire de $I$ avec les $P_j$ premiers. Alors nous avons
$$\prod_{k=1}^rx_k^{n_k}\in P_1\cap...\cap P_q+\m^{i+1}\ .$$
Donc pour tout $j$, $\prod_{k=1}^rx_k^{n_k}\in P_j+\m^{i+1}$. Donc d'apr\`es ce qui pr\'ec\`ede, il existe $k$ tel que $x_k\in P_j+\m^{\left\lfloor\frac{i-d}{c}\right\rfloor+1}$ avec $c$ et $d$ des constantes qui ne d\'ependent que des $P_j$. Fixons $k\in\{1,...,\,r\}$. Notons $J_k$ l'ensemble des $j$ tel que  $x_k\in P_j+\m^{\left\lfloor\frac{i-d}{c}\right\rfloor+1}$. Si $J_k=\emptyset$, nous posons alors $\ovl{x}_k=x_k$. Dans le cas contraire, pour tout $j\in J_k$,
$$x_k=\sum_{l\in H_j} p_{j,l}x_{j,l}+m_{k,j}$$
o\`u les  $p_{j,l}$ (quand $l$ parcourt l'ensemble $H_j$) engendrent $P_j$ et $m_{k,j}\in\m^{\left\lfloor\frac{i-d}{c}\right\rfloor+1}$ pour tout $j$. Soit $l_{j_1,j_2}$ la forme lin\'eaire
$$l_{j_1,j_2}(X_{j_1,l},\,X_{j_2,l'}):=\sum_{l\in H_{j_1}} p_{j_1,l}X_{j_1,l}-\sum_{l'\in H_{j_2}} p_{j_2,l'}X_{j_2,l'}\ .$$
Nous avons $l_{j_1,j_2}(x_{j_1,l},\,x_{j_2,l'})\in \m^{\left\lfloor\frac{i-d}{c}\right\rfloor+1}$ pour tout $j_1$ et $j_2$ dans $J_k$. D'apr\`es le th\'eor\`eme \ref{lin}, pour tous $j\in J_k$ et pour tout $l\in H_{j}$, il existe donc des $\ovl{x}_{j,l}\in x_{j,l}+\m^{\left\lfloor\frac{i-d'}{c'}\right\rfloor+1}$ tels que : 
$$l_{j_1,j_2}(\ovl{x}_{j_1,l},\ovl{x}_{j_2,l'})=0\text{ pour tout }j_1,\,j_2\in J_k, \text{ tout }l\in H_{j-1}\text{  et tout }l'\in H_{j_2},$$ avec $c'$ et $d'$ des constantes qui ne d\'ependent que des $P_j$. Nous notons alors $\ovl{x}_k=\sum p_{j_1,l}\ovl{x}_{j_1,l}$ et d'apr\`es ce qui pr\'ec\`ede 
$$\forall k\quad\ovl{x}_k\in \left(\bigcap_{j\in J_k}P_j\right)\cap (x_k+\m^{\left\lfloor\frac{i-d'}{c'}\right\rfloor+1})\ .$$
Comme $\cup_k J_k=\{1,...,\,r\}$, nous avons  
$$\prod_{k=1}^r\ovl{x}_k^{n_k}\in I\cap (\prod_{k=1}^rx_k^{n_k}+\m^{\left\lfloor\frac{i-d'}{c'}\right\rfloor+1})\ .$$
Donc il existe des $z^*_j$ tels que $\prod_{k=1}^r\ovl{x}_k^{n_k}+\sum_{j=1}^pf_jz^*_j=0$ ou encore
$$\sum_{j=1}^pf_j(z^*_j-z_j)\in\m^{\left\lfloor\frac{i-d'}{c'}\right\rfloor+1}\ .$$
Donc, d'apr\`es le lemme d'Artin-Rees, il existe des $\e_j\in\m^{\left\lfloor\frac{i-d''}{c''}\right\rfloor+1}$ tels que $\sum f_j(z^*_j-z_j+\e_j)=0$, o\`u $c''$ et $d''$ ne d\'ependent que des $P_j$ et de $I$. Nous posons alors $\ovl{z}_j=z_j-\e_j$ pour tout $j$. Nous avons donc 
$$\prod_{k=1}^r\ovl{x}_k^{n_k}+\sum_{j=1}^pf_j\ovl{z}_j=0$$
et $$\forall j\ \forall k,\ \ \ovl{x}_k-x_k,\, \ovl{z}_j-z_j\in\m^{\left\lfloor\frac{i-d''}{c''}\right\rfloor+1}.\qquad\Box$$

\begin{ex}
Soit $f$ un germe de fonction de Nash (resp. de fonction holomorphe). Alors si $f=gh$ avec $g$ et $h$ deux s\'eries formelles non inversibles alors $f$ peut s'\'ecrire comme le produit de deux germes de  fonctions de Nash (resp.  de deux fonctions holomorphes) non inversibles.\\
\end{ex}

\begin{ex}
Il est en g\'en\'eral faux que $XY$ admette une fonction de Artin. Consid\'erons par exemple l'anneau 
$$A:=\frac{\k[T_1,\,T_2]_{(T_1,T_2)}}{T_1^2-T_2^2(1+T_2)}\quad\text{ avec }\k\text{ un corps de caract\'eristique nulle}.\qquad$$
$A$ est irr\'eductible mais pas analytiquement irr\'eductible. Nous avons la relation $T_1^2-T_2^2(1+T_2)=(T_1-T_2\sqrt{1+T_2})(T_1+T_2\sqrt{1+T_2})$
o\`u $\sqrt{1+T_2}$ est une des deux s\'eries formelles dont le carr\'e vaut $1+T_2$. Soit $\left(\sqrt{1+T_2}\right)_n$ la s\'erie $\sqrt{1+T_2}$ tronqu\'ee \`a l'ordre $n$. Nous avons
$$\ord\left(\left(\sqrt{1+T_2}\right)_n-\sqrt{1+T_2}\right)=n+1.$$
Regardons le polyn\^ome $g(X,\,Y,\,Z)=XY-(T_1^2-T_2^2(1+T_2))Z$ de l'anneau $\k[T_1,\,T_2]_{(T_1,T_2)}[X,\,Y,\,Z]$. Posons
$$x_n=T_1\left(T_1-T_2\left(\sqrt{1+T_2}\right)_n\right),\,y_n=T_1+T_2\left(\sqrt{1+T_2}\right)_n\text{ et } z=T_1.$$
Nous avons $x_ny_n-(T_1^2-T_2^2(1+T_2))z\in\m^{n+4}$ pour tout entier $n\geq 1$. Or $x_n\notin (T_1^2-T_2^2(1+T_2))+\m^3$ et $y_n\notin (T_1^2-T_2^2(1+T_2))+\m^2$. Donc il n'existe pas de solution de $g$ ``proche'' de $(x_n,y_n,z)$ pour la topologie $\m$-adique.

\end{ex}
La preuve pr\'ec\'edente est constructive, dans le sens o\`u l'on peut donner une expression d'une fonction affine bornant la fonction de Artin de $g$ en terme de coefficients apparaissant dans des ICL  et de coefficients pour lesquels le lemme d'Artin-Rees est v\'erifi\'e pour des id\'eaux d\'ependants de $I$. N\'eanmoins ces bornes peuvent \^etre am\'elior\'ees \`a l'aide d'un th\'eor\`eme d\^u \`a D. Rees. Nous donnons un exemple ci-dessous.

\subsection{Bornes explicites}
Nous allons donner ici deux majorations affines de la fonction de Artin du polyn\^ome $X^n+\sum_jf_jZ_j$ : l'une \`a l'aide du th\'eor\`eme d'Izumi et l'autre \`a l'aide d'un th\'eor\`eme de Rees (cf. th\'eor\`eme \ref{Ree}).
\begin{lemme}
Soient $A$ un anneau local n{\oe}th\'erien complet et $I$ un id\'eal radical de $A$ engendr\'e par $f_1,...,\,f_p$. Soit $g(X,\,Z_j):=X^n+\sum_jf_jZ_j$.
Alors $g$ admet une fonction de Artin major\'ee par
$$i\lgm (2a)^{\left\lfloor\ln_2(n)\right\rfloor+1}(i+i_P+i_I)+b(1+2a+\cdots+(2a)^{\left\lfloor\ln_2(n)\right\rfloor})$$
o\`u $a$ et $b$ sont les plus petites constantes d'une ICL v\'erifi\'ee par tous les id\'eaux premiers associ\'es \`a $I$, $i_P$ est la plus petite constante pour laquelle le lemme d'Artin-Rees est v\'erifi\'e pour les id\'eaux engendr\'e par deux id\'eaux premiers associ\'es \`a $I$ et $i_I$ est la plus petite constante pour laquelle le lemme d'Artin-Rees est v\'erifi\'e pour $I$ (c'est-\`a-dire $I\cap \m^{i+i_I}\subset I\m^{i}$).\\
\end{lemme}

\textbf{Preuve :}
Soient $x$ et des $z_j$ tels que 
$$x^n+\sum_jf_jz_j\in\m^{(2a)^{\left\lfloor\ln_2(n)\right\rfloor+1}(i+i_P+i_I)+b(1+2a+\cdots+(2a)^{\left\lfloor\ln_2(n)\right\rfloor})+1}.$$ 
Soit $$I=P_1\cap\cdots\cap P_r$$
la d\'ecomposition primaire de $I$ avec les $P_l$ premiers.
Alors 
$$\nu_{P_l}(x^n)\geq (2a)^{\left\lfloor\ln_2(n)\right\rfloor+1}(i+i_P+i_I)+b(1+2a+\cdots+(2a)^{\left\lfloor\ln_2(n)\right\rfloor})+1$$
pour tout $l$.\\
Nous pouvons construire la suite suivante par r\'ecurrence (o\`u $n_0=n$) :\\
Si $n_k$ est pair on pose $n_{k+1}=\frac{n_k}{2}$, sinon on pose $n_{k+1}=\frac{n_k+1}{2}$.
Ecrivons $n_k$ et $n_{k+1}$ en base 2 :
$$n_k=\a_0+\a_12+\cdots+\a_{q-1}2^{q-1}+2^q\quad(q=\left\lfloor\ln_2(n_k)\right\rfloor)$$
$$n_{k+1}=\b_0+\b_12+\cdots+\b_{q-1}2^{q-1}+\b_q2^q$$
avec les $\a_j$ et les $\b_j$ dans $\{0,1\}$.\\
Si $\a_0=0$, alors $\b_q=0$ et $\b_{q-1}=1$. Si $\a_0=\a_1=\cdots=\a_{q-1}=1$ alors $\b_0=\b_1=\cdots=\b_{q-1}=0$ et $\b_q=1$. Si l'un des $\a_j$, pour $0\leq j\leq q-1$, est nul, alors $\b_q=0$.\\
Si $\a_0=\a_1=\cdots=\a_{q-1}=0$  alors $\b_0=\b_1=\cdots=\b_{q-2}=0$ et $\b_{q-1}=1$. Nous voyons donc, si $q=\left\lfloor\ln_2(n)\right\rfloor$, que $n_{q}=1\text{ ou } n_{q+1}=1$.\\
Donc, d'apr\`es les hypoth\`eses, nous avons
$$\nu_{P_l}(x^{n_1})\geq (2a)^{\left\lfloor\ln_2(n)\right\rfloor}(i+i_P+i_I)+b(1+2a+\cdots+(2a)^{\left\lfloor\ln_2(n)\right\rfloor-1})+1\ .$$
Par induction nous avons alors 
$$\nu_{P_l}(x)\geq i+i_P+i_I+1\ .$$ 
Il existe donc des $x_{l,j}$ tels que $x-\sum_jp_{l,j}x_{l,j}\in\m^{i+i_P+i_I+1}$ o\`u les $p_{l,j}$ engendrent $P_l$. D'apr\`es la derni\`ere \'etape de la preuve du th\'eor\`eme \ref{thmprinc}, il existe donc $\ovl{x}\in (P_1\cap...\cap P_r)\cap\left( x+\m^{i+i_I+1}\right)$.\\
Il existe alors des $z^*_j$ tels que $\ovl{x}=\sum_jf_jz^*_j$ et $x-\sum_jf_jz^*_j\in\m^{i+i_I+1}$. Notons  $\ovl{x}^n=\sum_jf_jz^{**}_j$ avec les $z^{**}_j$ dans $A$. Nous avons alors $\sum_jf_j(z_j+z^{**}_j)\in\m^{i+i_I+1}$ et il existe alors des $t_j\in z_j+z^{**}_j+\m^{i+1}$ tels que $\sum_jf_jt_j=0$. On pose alors $\ovl{z}_j=t_j-z^{**}_j$. Nous avons bien $g(\ovl{x},\,\ovl{z}_j)=0$ et $\ovl{x}-x\in\m^{i+1}$ et $\ovl{z}_j-z_j\in\m^{i+1}$ pour tout $j$.
$\quad\Box$\\
\\
Nous voyons ici que le coefficient $\l$ de la fonction $i\lgw \l i+c$ d\'ecrite ci-dessus est de la forme $n^c$ pour une constante $c\geq 1$. Il est possible dans ce cas d'am\'eliorer cette borne \`a l'aide du th\'eor\`eme suivant :\\
\begin{theoreme}\label{Ree}\cite{Re3}
Soit $A$ un anneau local et n{\oe}th\'erien et $I$ un id\'eal de $A$ tel que $A/I$ est non ramifi\'e. Alors, pour tout $x$ dans $A$, la limite $\lim_n\frac{\nu_I(x^n)}{n}$ existe et est \'egale \`a la limite sup\'erieure de cette suite. Notons $\ovl{\nu}_I$ la fonction d\'efinie par 
$$\forall x\in A, \ \ovl{\nu}_I(x)=\lim_n\frac{\nu_I(x^n)}{n}.$$
Il existe alors une constante $c\geq  0$ telle que
$$\forall x\in A, \ \nu_I(x)\leq  \ovl{\nu}_I(x)\leq \nu_I(x)+c.$$
\end{theoreme}
Pour un entier $c$ nous notons $\lceil c\rceil$ sa partie enti\`ere sup\'erieure, c'est-\`a-dire $\lceil c\rceil=c$ si $c$ est entier et $\lceil c\rceil=\left\lfloor c\right\rfloor+1$ si $c$ n'est pas entier. Nous pouvons alors d\'eduire le lemme suivant
\begin{lemme}\label{xn}
Soit $A$ un anneau local n{\oe}th\'erien complet et $I$ un id\'eal radical de $A$ engendr\'e par $f_1,...,\,f_p$. Soit $g(X,\,Z_j):=X^n+\sum_jf_jZ_j$.
Alors $g$ admet une fonction de Artin major\'ee par la fonction  
$$i\lgm n\left\lceil \frac{i+i_I}{n}\right\rceil +nc\leq i+i_I+n(c+1)$$
 o\`u $c$ est la plus petite constante telle que $\forall x\in A, \ \ovl{\nu}_I(x)\leq \nu_I(x)+c$ et $i_I$ est la plus petite constante pour laquelle le lemme d'Artin-Rees est v\'erifi\'e pour $I$ (c'est-\`a-dire $I\cap \m^{i+i_I}\subset I\m^{i}$).\\

\end{lemme}
\textbf{Preuve :}
Soient $x$ et des $z_j$ tels que 
$$x^n+\sum_jf_jz_j\in\m^{n\left\lceil \frac{i+i_I}{n}\right\rceil +nc+1}$$
avec les notations du lemme. Alors 
$$\frac{\nu_I(x^n)}{n}\leq \ovl{\nu}_I(x) \leq \nu_I(x)+c$$ d'apr\`es le th\'eor\`eme de Rees. Or nous avons  $\nu_I(x^n)\geq n\left\lceil \frac{i+i_I}{n}\right\rceil +nc+1$, donc $\nu_I(x)\geq \left\lceil \frac{i+i_I}{n}\right\rceil+1$. Il existe alors des $z^*_j$ tels que $x-\sum_jf_jz^*_j\in\m^{\left\lceil \frac{i+i_I}{n}\right\rceil+1}$, c'est-\`a-dire $x=\sum_jf_jz^*_j+\e$ avec $\e\in\m^{\left\lceil \frac{i+i_I}{n}\right\rceil+1}$.  D'o\`u $x^n=\sum_jf_jR_j(z^*_j,\,\e)+\e^n$ avec $R_j$ des polyn\^omes en $p+1$ variables. D'o\`u
$\sum_jf_j(z_j+R_j(z^*_j,\,\e))\in\m^{i+i_I+1}$ et il existe alors des $t_j\in z_j+R_j(z^*_j,\,\e)+\m^{i+1}$ tels que $\sum_jf_jt_j=0$. On pose alors $\ovl{z}_j=t_j-R_j(z^*_j,\,\e)$ et $\ovl{x}=\sum_jf_jz^*_j$.
$\quad\Box$\\
\\
Nous allons maintenant utiliser ce dernier lemme pour obtenir des d\'eterminations explicites de cl\^otures int\'egrales approch\'ees d'id\'eaux.

\section{Application \`a des d\'eterminations explicites de cl\^otures int\'egrales approch\'ees d'id\'eaux}\label{clotint}
\subsection{Cl\^oture int\'egrale approch\'ee d'un id\'eal }
Nous commen\c{c}ons tout d'abord par rappeler certains r\'esultats connus. Si $I$  est un id\'eal d'un anneau $A$ int\`egre, nous notons $\ovl{I}$ sa cl\^oture int\'egrale. Il est bien connu (voir par exemple \cite{Eis}) que $I\subset\ovl{I}\subset\sqrt{I}$. En particulier si $I$ est radical alors $\ovl{I}=I$. D'autre part si $A$ est principal, alors $\ovl{I}=I$.\\
D. Delfino et I. Swanson ont montr\'e le th\'eor\`eme suivant qui est une g\'en\'eralisation d'un th\'eor\`eme de Rees \cite{Re1}:
\begin{theoreme}\label{D-S}\cite{D-S}
Soit $(A,\m)$ un anneau local n{\oe}th\'erien excellent. Soit $I$ un id\'eal de A. Alors il existe $a$ et $b$ des entiers tels que 
$$\ovl{I+\m^{ai+b}}\subset \ovl{I}+\m^i\qquad\forall i\in\N$$
$$\text{ou  encore }\ \ovl{I+\m^i}\subset \ovl{I}+\m^{\left\lfloor\frac{i-b}{a}\right\rfloor}\ \ \forall i\in\N.$$
\end{theoreme}
Pour prouver ce th\'eor\`eme, D. Delfino et I. Swanson se ram\`enent au cas o\`u $I$ est principal et $A$ complet et normal. Dans ce cas elles montrent que 
 tout \'el\'ement de $\ovl{I+\m^{i}}$ v\'erifie une relation de la forme
$$X^n+X^{n-1}\sum_jg_jX_{1,j}+\cdots+\sum_{j_1\leq\cdots\leq j_n}g_{j_1}...g_{j_n}X_{n,j_1,...,j_n}\in \m^{\left\lfloor\frac{i}{l}\right\rfloor}$$
o\`u $n$ et $l$ sont ind\'ependants de l'\'el\'ement choisi et de l'entier $i$. Ensuite elles montrent, toujours sous les m\^emes hypoth\`eses ($I$ principal et $A$ complet et normal), que le polyn\^ome pr\'ec\'edent admet une fonction de Artin major\'ee par une fonction affine (th\'eor\`eme 3.10 de \cite{D-S}).\\
Nous allons donner dans cette partie  une g\'en\'eralisation du th\'eor\`eme 3.10 de \cite{D-S}. L'int\'er\^et de notre preuve vient du fait que celle-ci est constructive et permet d'obtenir des bornes explicites  en termes de coefficients apparaissant dans certaines ICL.
\subsection{G\'en\'eralisation d'un r\'esultat de Delfino et Swanson}
En utilisant le lemme \ref{xn}, nous allons donc donner deux propositions qui g\'en\'eralisent le th\'eor\`eme 3.10 de \cite{D-S} :
\begin{proposition}\label{DS1}
Soit $$g(X,\,X_{1,j},...,\,X_{n,j_1,...,j_n},\,Y_1,...,\,Y_q):=X^n+X^{n-1}\sum_jg_jX_{1,j}+\cdots$$
$$+\sum_{j_1\leq\cdots\leq j_n}g_{j_1}...g_{j_n}X_{n,j_1,...,j_n}+\sum_{l=1}^qf_lY_l\ $$
avec les $g_j$  et les $f_l$ dans $A$, local complet n{\oe}th\'erien, tels que $I=(f_l)+(g_j)$ soit radical.
Alors $g$ admet une fonction de Artin major\'ee par la fonction 
$$i\lgm i+i_I+n(c+1)$$ o\`u $c$ est la plus petite constante telle que $\forall x\in A, \ \ovl{\nu}_I(x)\leq \nu_I(x)+c$, et $i_I$ est la plus petite constante pour laquelle le lemme d'Artin-Rees est v\'erifi\'e pour $I$ (c'est-\`a-dire $I\cap \m^{i+i_I}\subset I\m^{i}$).\\
\end{proposition}

\textbf{Preuve :}
Soient $(x,\,x_{1,j},...,\,x_{n,j_1,...,j_n},\,y_1,...,\,y_q)\in A$ tels que 
$$g(x,\,x_{1,j},...,\,x_{n,j_1,...,j_n},\,y_l)\in\m^{i+i_I+n(c+1) }\ .$$ 
Posons 
$$t'_j=x^{n-1}x_{1,j}+x^{n-2}\sum_{j_2\geq j}g_{j_2}x_{2,j,j_2}+\cdots+\sum_{j_n\geq \cdots\geq j_2\geq j}g_{j_2}...g_{j_n}x_{n,j,j_2,...,j_n}\ .$$
Alors nous avons 
$$x^n+\sum_jg_jt'_j+\sum_lf_ly_l=g(x,\,x_{1,j},...,\,x_{n,j_1,...,j_n},\,y)\ .$$
D'apr\`es la preuve du  lemme \ref{xn}, il existe $\ovl{x}\in x+\m^{i+i_{I}+1}$ tel que $\ovl{x}\in I$. Nous pouvons \'ecrire $\ovl{x}=\sum_jg_jx'_j+\sum_l f_lz_l$. Nous avons alors
$$g(\ovl{x},\,x_{1,j},...,\,x_{n,j_1,...,j_n},\,y_l)\in\m^{i+i_{I}+1}.$$
D'o\`u
$$\sum_{j_1\leq\cdots\leq j_n}g_{j_1}...g_{j_n}\left(x_{n,j_1,...,j_n}+h_{j_1,...,j_n}(x_{1,j},...,x_{n-1,j'_1,...,j'_{n-1}},x'_j)\right)+\qquad\qquad\qquad$$
$$\qquad\qquad\qquad\qquad\qquad\qquad\qquad\qquad+\sum_l f_lt_l\in\m^{i+i_{I}+1}$$
avec $t_l=y_l+t^*_l(x'_j,z_l)$ et $h_{j_1,...,j_n}$ polynomiale \`a coefficients dans $A$. D'apr\`es Artin-Rees, il existe alors $(\ovl{t}_1,...,\,\ovl{t}_q)\in (t_1,...,t_q)+\m^{i+1}$ et 
$$\ovl{t}_{j_1,...,j_n}\in x_{n,j_1,...,j_n}+h_{j_1,...,j_n}(x_{1,j},...,x_{n-1,j'_1,...,j'_{n-1}},x'_j) +\m^{i+1}$$
tels que $\sum_{j_1\leq\cdots\leq j_n}g_{j_1}...g_{j_n}\ovl{t}_{j_1,...,j_n}+\sum_l f_l\ovl{t}_l=0$.
Posons alors 
$$\ovl{x}_{i,j_1,...,j_i}=x_{i,j_1,...,j_i} \text{ pour tout }i<n$$
 et 
$$\ovl{x}_{n,j_1,...,j_n}=\ovl{t}_{j_1,...,j_n}-h_{j_1,...,j_n}(x_{1,j},...,x_{n-1,j'_1,...,j'_{n-1}},x'_j)\ .$$ 
Nous avons $\ovl{x}_{i,j_1,...,j_i}-x_{i,j_1,...,j_i}\in\m^{i+1}$ pour tout $i$ et $j_k$. Posons $\ovl{y}_l=\ovl{t}_l-t^*_l$ pour tout $l$. Nous avons donc $\ovl{y}_l-y_l\in\m^{i+1}$ et $\ovl{x}-x\in\m^{i+1}$. De plus il est clair que $g(\ovl{x},\ovl{x}_j,\ovl{y}_l)=0$.$\quad\Box$

\begin{proposition}\label{xnf}
Soit $$g(X,\,X_{1},...,\,X_{n},\,Y_1,...,\,Y_q)=X^n+f^tX^{n-1}X_1+\cdots+f^{nt}X_n+\sum_{l=1}^qf_lY_l\ $$
avec les  $f_j$ et $f$ dans $A$, local complet n{\oe}th\'erien, tels que $((f_j):f)=(f_j)$ et $I:=(f,\,f_j)$ soit radical, et soit $t$ un entier strictement positif. Alors $g$ admet une fonction de Artin major\'ee par
$$i\lgw i+i_I+ti_{J_n} +tn(c+1)$$ o\`u $c$ est la plus petite constante telle que $\forall x\in A, \ \ovl{\nu}_I(x)\leq \nu_I(x)+c$, et $i_{J_n}$ est la plus petite constante pour laquelle le lemme d'Artin-Rees est v\'erifi\'e pour $J_n=(f^n, (f_j))$ (c'est-\`a-dire $J_n\cap \m^{i+i_{J_n}}\subset J_n\m^{i}$) et $i_I$ est la plus petite constante pour laquelle le lemme d'Artin-Rees est v\'erifi\'e pour $I$.
\end{proposition}

\textbf{Preuve :} Notons $I:=(f,\,f_l)$.
Soit $i$ un entier positif. Soient $x$, des $x_j$ et des $y_k$ tels que 
$$g(x,\,x_j,\,y_k)\in\m^{i+i_I+ti_{J_n} +tn(c+1)+1}.$$ 
Supposons tout d'abord qu'il existe $\ovl{x}$ et des $\ovl{x}_j$ tels que $\ovl{x}-x$, $\ovl{x}_j-x_j\in\m^{i+i_I+1}$ pour tout $j$, et $\ovl{x}^n+f^t\ovl{x}^{n-1}\ovl{x}_1+\cdots+f^{nt}\ovl{x}_n\in (f_l)$. Il existe alors des $t_l$ tels que $\ovl{x}^n+f^t\ovl{x}^{n-1}\ovl{x}_1+\cdots+f^{nt}\ovl{x}_n=\sum_lf_lt_l$. Nous avons alors $\sum_lf_l(y_l+t_l)\in\m^{i+i_I+1}$. Il existe alors des $z_l\in y_l+t_l+\m^{i+1}$ tels que $\sum_lf_lz_l=0$. Nous posons alors $\ovl{y}_l=z_l-t_l$. Nous avons $\ovl{y}_l-y_l\in\m^{i+1}$ pour tout $l$ et $g(\ovl{x},\,\ovl{x}_j,\,\ovl{y}_l)=0$. Nous pouvons donc supposer que $(f_l)=(0)$.\\
Supposons alors $(f_l)=(0)$. Alors, comme dans la preuve de la proposition pr\'ec\'edente, il existe $\ovl{x}\in I$ tel que nous ayons $\ovl{x}=x$ modulo $\m^{i+ti_{J_n}+n(t-1)(c+1)+1}$. Donc nous avons
$$x=fx'+\e_1$$
avec $\e_1\in\m^{i+i_I+ti_{J_n}+n(t-1)(c+1)+1}$. Nous avons alors
$$f^nx'^n+f^{t+n-1}x'^{n-1}x_1+\cdots+f^{nt}x_n\in\m^{i+i_I+ti_{J_n}+n(t-1)(c+1)+1}.$$
 D'o\`u
$$f^n\left(x'^n+f^{t-1}x'^{n-1}x_1+\cdots+f^{n(t-1)}x_n\right)\in\m^{i+i_I+ti_{J_n}+n(t-1)(c+1)+1}\ .$$
Donc nous avons $x'^n+f^{t-1}x'^{n-1}x_1+\cdots+f^{n(t-1)}x_n\in\m^{i+i_I+(t-1)i_{J_n}+n(t-1)(c+1)+1}$, car $((0):f)=(0)$.

On obtient le r\'esultat par r\'ecurrence sur $t$, car pour $t=0$ le polyn\^ome est lisse en tout point (le coefficient de $x_n$ est \'egal \`a 1). $\quad\Box$

\subsection{Exemple effectif}\label{ex}
Cet exemple est cit\'e dans \cite{D-S} mais incorrectement \'etudi\'e car les auteurs utilisent un r\'esultat de M. Lejeune-Jalabert uniquement valable pour $A=\k[[T]]$. Pour \'etudier cet exemple, nous allons utiliser ici la proposition \ref{xnf} et un r\'esultat de Delfino et Swanson \cite{D-S}.\\
Soient $a,\,t,\,N\in\N$ tels que $a\geq 2$, $t\geq 1$ et $N\geq 3$ et $\k$ un corps contenant les racines $a$-i\`emes de l'unit\'e et de caract\'eristique ne divisant pas $a$. Notons 
$$A:=\frac{\k[[T_1,...,\,T_N]]}{(T_1^a+\cdots+T_N^a)}\ .$$
Soit $B=\k[[T_1,\,T_2,...,\,T_{N-1}]]$. L'extension $\Frac(A)\subset\Frac(B) $ est  galoisienne et s\'eparable et notons $n=[\Frac(A):\Frac(B)]$. L'entier $n$ divise $\Phi(a)$, la fonction d'Euler de $a$, donc $n<a$. Nous utilisons alors le
\begin{lemme}\cite{D-S}
Soit $(A,\m)$ un anneau local complet normal n{\oe}th\'erien et soit $f$ un \'el\'ement non nul de $A$. Soit $B=\k[[f,f_2,...,\,f_N]]$ o\`u $(f,f_2,...,\,f_N)$ est un syst\`eme de param\`etres de $A$. Supposons que $\Frac(A)\subset\Frac(B)$  est une extension galoisienne s\'eparable et notons $n=[\Frac(A):\Frac(B)]$. Alors tout \'el\'ement de $\ovl{f^tA+\m^i}$ v\'erifie une \'equation de degr\'e $n$ sur $f^tA+\m^{\left\lfloor\frac{i}{nt}\right\rfloor}$
\end{lemme}

Donc d'apr\`es le lemme pr\'ec\'edent, tout \'el\'ement de $\ovl{T_1^tA+\m^i}$ v\'erifie une \'equation de degr\'e $n$ sur $T_1^tA+\m^{\left\lfloor\frac{i}{nt}\right\rfloor}$.\\
Soit $x\in A$ v\'erifiant une \'equation de degr\'e $n$ sur $T_1^tA+\m^{\left\lfloor\frac{i}{nt}\right\rfloor}$. Notons $I$ l'id\'eal $(T_1,\,T_2^a+\cdots+T_N^a)$. Si $N>3$, $\nu_I$ est une valuation  car $Gr_{\m}A/I$ est int\`egre.\\
Si $N\geq 3$, l'id\'eal $I$ \'etant homog\`ene et radical, nous avons aussi $c=0$. D'apr\`es le corollaire \ref{coefinit}, nous avons $i_{J_n}=a$.\\
Si $N\geq 3$, d'apr\`es le lemme \ref{lininiirr}, comme $T_1$ et $T_2^a+\cdots+T_N$ forment une suite r\'eguli\`ere, $i_I=a$.\\
Donc, d'apr\`es la proposition \ref{xnf}, il existe 
$\ovl{x}\in T_1^tA\cap\left(x+\m^{\left\lfloor\frac{i-a}{nt}\right\rfloor-t(a+n)}\right).$
Nous obtenons alors la

\begin{proposition}
Soient $a,\,t,\,N\in\N$ tels que $a\geq 2$, $t\geq 1$ et $N\geq 3$ et $\k$ un corps contenant les racines $a$-i\`emes de l'unit\'e et de caract\'eristique ne divisant pas $a$ et $A=\frac{\k[[T_1,...,\,T_N]]}{(T_1^a+\cdots+T_N^a)}$. Alors
\begin{equation}
\forall i\in\N^* \ \ \qquad\ovl{T_1^tA+\m^i}\subset T_1^tA+\m^{\left\lfloor\frac{i-a}{nt}\right\rfloor-t(a+n)}
\end{equation}
o\`u $n=[\Frac(A):\Frac(B)]$.\\
\\
\end{proposition}

\end{document}